\documentclass[11pt]{article}
\usepackage[margin =1in]{geometry}
\usepackage{amssymb,amsthm,amsmath}
\usepackage{enumerate}
\usepackage{hyperref}
\usepackage{cite}
\usepackage[capitalize]{cleveref}
\usepackage{enumitem}
\usepackage{color}
\usepackage{cancel}
\usepackage{pgf}
\usepackage{svg}
\usepackage{pgfplots}
\usepackage{import}
\usepackage{float}
\usepackage{ulem}

\usepackage[utf8]{inputenc} 
\usepackage[T1]{fontenc}    
\usepackage{hyperref}       
\usepackage{url}            
\usepackage{booktabs}       
\usepackage{amsfonts}       
\usepackage{nicefrac}       
\usepackage{microtype}
\usepackage{multirow}
\usepackage{xfrac}
\usepackage{todonotes}
\usepackage{titlesec}
\usepackage{caption}
\usepackage{subcaption}
\usepackage{algorithm}
\usepackage{algorithmic}

\titleformat{\subsubsection}[runin]
{\normalfont\normalsize\bfseries}{\thesubsubsection}{1em}{}

\usepackage{graphicx}

\usepackage{appendix}
\numberwithin{equation}{section}

\newcommand{\inclu}[0] {\ar@{^{(}->}}

\newcommand{\dist}{{\rm dist}}

\newcommand{\prox}{\text{prox}}

\newcommand{\proj}{\mathrm{proj}}

\newcommand{\argmin}{\operatornamewithlimits{argmin}}

\newtheorem{thm}{Theorem}[section]
\newtheorem{theorem}{Theorem}[section]
\newtheorem{definition}[thm]{Definition}

\newtheorem{lemma}[thm]{Lemma}

\newtheorem{corollary}[thm]{Corollary}

\newtheorem{assumption}{Assumption}

\crefname{claim}{claim}{claims}
\Crefname{claim}{Claim}{Claims}
\crefname{lem}{lemma}{lemmas}
\Crefname{lem}{Lemma}{Lemmas}
\crefname{algorithm}{algorithm}{algorithms}
\Crefname{algorithm}{Algorithm}{Algorithms}

\newtheorem{remark}{Remark}

\theoremstyle{remark}

\newcommand{\blk}{black}

\usepackage{mathtools}

\usepackage[algo2e, boxruled]{algorithm2e}

\usepackage[T1]{fontenc}
\usepackage{lmodern}
\allowdisplaybreaks
\pgfplotsset{compat=1.14}

\begin{document}

	\title{First-Order Methods for Nonsmooth Nonconvex Functional Constrained Optimization with or without Slater Points }

	 \author{Zhichao Jia\footnote{Johns Hopkins University, Department of Applied Mathematics and Statistics, \url{zjia12@jhu.edu}} \qquad Benjamin Grimmer\footnote{Johns Hopkins University, Department of Applied Mathematics and Statistics, \url{grimmer@jhu.edu}}}

	\date{}
	\maketitle

	\begin{abstract}
	    Constrained optimization problems where both the objective and constraints may be nonsmooth and nonconvex arise across many learning and data science settings. In this paper, we show for any Lipschitz, weakly convex objectives and constraints, a simple first-order method finds a feasible, $\epsilon$-stationary point at a convergence rate of $O(\epsilon^{-4})$ without relying on compactness or Constraint Qualification (CQ). When CQ holds, this convergence is measured by approximately satisfying the Karush–Kuhn–Tucker conditions. When CQ fails, we guarantee the attainment of weaker Fritz-John conditions. As an illustrative example, our method stably converges on piecewise quadratic SCAD regularized problems despite frequent violations of constraint qualification. The considered algorithm is similar to those of~\cite{ma2020quadratically,boob2022stochastic} (whose guarantees further assume compactness and CQ), iteratively taking inexact proximal steps, computed via an inner loop applying a switching subgradient method to a strongly convex constrained subproblem. Our non-Lipschitz analysis of the switching subgradient method appears to be new and may be of independent interest.
	    
    \end{abstract} 
    
    \section{Introduction}

In this paper, we considered the difficult family of constrained optimization problems where both the objective and constraints may be nonconvex and nonsmooth. Specifically, we consider problems of the following form:
\begin{equation}
\begin{cases}
    \min_{x \in X} \quad &f(x) \\
    \mathrm{s.t.} \quad &g_i(x) \leq 0, \qquad i=1,...,m.
    \label{mainproblemoriginal}
\end{cases}
\end{equation}
for some closed convex domain $X \subseteq \mathbb{R}^d$. The objective $f:X \to \mathbb{R}$ and constraints $g_i:X \to \mathbb{R}, i=1,...,m$ are assumed to be continuous on $X$ {\color{black} and weakly convex (see Section~\ref{preliminary})}, but need not be convex nor differentiable.

Constrained optimization problems with nonsmooth and nonconvex objectives and constraints are common in modern data science and machine learning. For instance, phase retrieval, blind deconvolution, and covariance matrix estimation all fall within nonconvex and nonsmooth minimization~\cite{davis2020nonsmooth,duchi2019solving,charisopoulos2021composite,chen2015exact,davis2019stochastic}. If sparsity of solutions is expected or desired, often a regularizing constraint is introduced (e.g., convex choices like $\ell_1$-norms or $\ell_2$-norms, nonconvex choices like SCAD functions~\cite{antoniadis1997wavelets,fan2001variable} or $\ell_q$-norms for $q\in(0,1)$). SCAD functions will serve as a running example throughout this work as they are simple piecewise quadratic functions exhibiting nonsmoothness and nonconvexity, with widespread usage~\cite{kim2008smoothly,fan2009network,xie2009scad,gasso2009recovering,zeng2022moreau}. Other problems like multi-class Neyman-Pearson classification~\cite{weston1998multi,ma2020quadratically,tian2021neyman}, minimizing the loss on one class while controlling the losses on other classes under some values, provide another typical setting of constrained optimization inheriting any nonsmoothness and nonconvexities from the loss functions.

Our approach to solving nonsmooth, nonconvex, constrained problems relies on two classic ingredients, which we review below: (inexact) proximal point methods and Fritz-John/Karush-Kuhn-Tucker stationarity conditions.

\paragraph{(Inexact) Proximal Point Methods} {\color{black} Based on the concept of the proximal mapping introduced in~\cite{moreau1965proximite}, proximal point methods, explored initially by~\cite{martinet1970regularisation,guler1992new}, repeatedly minimize a quadratically penalized objective. Namely, supposing $m=0$ (i.e., no  functional constraints are given), one can seek a stationary point of $\min_{x\in X} f(x)$ by iterating 
\begin{equation}
\begin{aligned}
    x_{k+1} = \prox_{\alpha, f}(x_k):= \argmin_{x \in X}\left\{ f(x)+\frac{1}{2\alpha}\|x-x_k\|^2\right\} 
    \label{PPM}
\end{aligned}
\end{equation}
with stepsize $\alpha>0$.
Many works~\cite{hare2009computing,hare2010redistributed,salzo2012inexact,paquette2018catalyst,davis2019stochastic,davis2019proximally,bertsekas2015convex,frostig2015regularizing,kong2019complexity,lan2019accelerated} since have continued this line of work, solving {\color{black} weakly convex} problems via repeated inexact evaluation of the above proximal operator.}
By restricting to the family of weakly convex functions (defined in~\eqref{weaklyconvex}), this proximal subproblem is guaranteed to be convex with a unique solution for small enough $\alpha$. When the proximal map can be evaluated exactly, an $\epsilon>0$-stationary point (defined in Definitions \ref{definitionFJ1} and \ref{definitionKKT1}) is found within $O(1/\epsilon^2)$ iterations. The inexact methods of~\cite{davis2019proximally,davis2019stochastic} show that using cheaper subgradient oracle calls such a point is found within $O(1/\epsilon^4)$ iterations.

We follow the extension of these ideas to nonconvex ({\color{black} weakly convex}) inequality constraints proposed by Ma et al.~\cite{ma2020quadratically} and Boob et al.~\cite{boob2022stochastic}. Their ideas and comparisons with our contributions are discussed in Section~\ref{relatedwork}. To this end, we consider the following proximal subproblem, penalizing the constraints in addition to the objective
\begin{equation}
\begin{aligned}
    x_{k+1} \approx \argmin_{x \in X}\left\{ f(x)+\frac{1}{2\alpha}\|x-x_k\|^2 \mid g_i(x)+\frac{1}{2\alpha}\|x-x_k\|^2 \leq \tau {\color{\blk} \ \forall i=1,...,m} \right\}
    \label{conssubproblem}
\end{aligned}
\end{equation}
with stepsize $\alpha>0$ and feasibility tolerance $\tau\geq 0$.
Importantly, any feasible solution to this proximal subproblem $x_{k+1}$ has its infeasibility bounded by
$ g_i(x_{k+1}) \leq \tau - \frac{1}{2\alpha}\|x_k-x_{k+1}\|^2 $.
Hence a sequence of $x_k$ generated by inexactly evaluating this mapping remains feasible for the original problem~\eqref{mainproblemoriginal} until it reaches approximate stationarity (that is, $\|x_k-x_{k+1}\| \geq \sqrt{2\alpha\tau}$ implies $g_i(x_{k+1})\leq 0$ for each constraint $i$).

\paragraph{Fritz-John/Karush-Kuhn-Tucker Stationarity} Let $\partial f(x)$ denote a generalized subdifferential of a function $f$ and $N_X(x)$ denote the normal cone of $X$ at $x$, formally defined in Section~\ref{preliminary}. Here, we consider two classic measurements of stationarity: Fritz-John (FJ) conditions giving a weaker optimality condition and Karush-Kuhn-Tucker (KKT) conditions giving a stronger condition. {\color{black} These classic optimality conditions are formalized below without reliance on differentiability of $f$ or $g_i$; see~\cite[Section 3.1]{Clark} for such a development more generally.}

We say that a feasible solution $x^*$ is a FJ point of~\eqref{mainproblemoriginal} if there exists nonnegative multipliers $\gamma_0^* \in \mathbb{R}$ and $\gamma^*=(\gamma_1^*,...,\gamma_m^*)^T \in \mathbb{R}^m$, and subgradients $\zeta_f \in \partial f(x^*)$ and $\zeta_{gi} \in \partial g_i(x^*)$ such that $(\gamma_0^*,\gamma_1^*,...,\gamma_m^*)$ is a non-zero vector with
\begin{equation}
\begin{aligned}
    \gamma_i^*g_i(x^*)&=0, \qquad \forall i=1,...,m, \\
    \gamma_0^*\zeta_f+\sum_{i=1}^m\gamma_i^*\zeta_{gi} &\in -N_X(x^*).
    \label{FJoriginal}
\end{aligned}
\end{equation}

Note requiring $(\gamma_0^*, \gamma_1^*, ..., \gamma_m^*)$ to be a nonzero vector could be equivalently expressed as requiring $\gamma_0^*+\sum_{i=1}^m\gamma_i^*=1$. This condition is necessary for $x^*$ to be a global (or local) minimizer~\cite{john2014extremum}. However, this condition can only give limited insights into the quality of $x^*$ as a solution when $\gamma_0^*=0$ since~\eqref{FJoriginal} becomes independent of $f$~\cite{mangasarian1967fritz}. {\color{black} In such a situation, FJ corresponds to being a stationary point of the function $\max\{g_1(x)\dots g_m(x)\}$ with value zero. That is, up to first-order, no direction strictly improves the feasibility of $x^*$. Although not a desirable position, this is a sensible place for a first-order method to halt since it cannot identify any feasible direction to move. For example, if $x_k$ is a FJ point, the proximal step~\eqref{conssubproblem} with $\tau=0$ would remain stationary $x_k$ is the only feasible solution.}

This weakness is remedied by the stronger notion of KKT points, which implicitly require $\gamma_0^*\neq 0$. We say a feasible $x^*$ is a KKT point for the problem~\eqref{mainproblemoriginal} if there exists nonnegative Lagrange multipliers $\lambda^*{\color{\blk}=(\lambda_1^*,...,\lambda_m^*)^T}\in\mathbb{R}^m$, $\zeta_f \in \partial f(x^*)$ and $\zeta_{gi} \in \partial g_i(x^*)$ such that
\begin{equation}
\begin{aligned}
    \lambda_i^*g_i(x^*)&=0, \qquad \forall i=1,...,m, \\
    \zeta_f+\sum_{i=1}^m\lambda_i^*\zeta_{gi} &\in -N_X(x^*).
    \label{KKToriginal}
\end{aligned}
\end{equation}

The KKT conditions strengthen FJ, requiring $\gamma_0^*\neq 0$, in particular $\gamma_0^*=1$. The requirement that $\gamma_0^*\neq 0$ is equivalent to having the Mangasarian-Fromovitz Constraint Qualification (MFCQ) condition {\color{black} (see~\cite{MANGASARIAN196737} for smooth optimization)} hold: {\color{black} Denote the dual of the normal cone by $N^*_X(x)$, and} let $A(x)=\{i \mid g_i(x) = 0, i=1,...,m\}$. We say MFCQ holds at $x^*$ if
\begin{equation}
\begin{aligned}
    \exists v \in -N_X^*(x^*) \qquad s.t. \quad \zeta_{gi}^Tv<0 \quad \forall i \in A(x^*), \forall \zeta_{gi} \in \partial g_i(x^*).
    \label{MFCQ}
\end{aligned}
\end{equation}
{\color{black} Constraint qualifications such as MFCQ applied in nonsmooth problems were discussed in~\cite{jourani1994constraint} based on the Clarke subdifferential and can also be found in~\cite{boob2022stochastic}.} Approximate FJ and KKT stationarity measurements can differ greatly when constraint qualification does not hold. When a strengthened ($\sigma$-strong) MFCQ condition (defined later as~\eqref{strongMFCQ}) is satisfied, we can uniformly bound the size of any associated Lagrange multipliers. Without this, these multipliers may be arbitrarily large, even failing to exist when MFCQ fails. Consequently, approximate KKT stationarity may never be attained despite the iterates $x_k$ of~\eqref{conssubproblem} converging. 
In contrast, we show FJ conditions are approximately satisfied whenever $x_k$ converges.

{\color{black} This distinction motivates the development of convergence theory capable of guaranteeing FJ conditions in general and KKT conditions when constraint qualification holds. Prior works for smooth constrained optimization have developed such two-pronged theory~\cite{cartis2011evaluation,cartis2013evaluation,cartis2014complexity,birgin2016evaluation,hinder2018worst}, ensuring their methods produce approximate {\color{\blk} scaled} KKT points (a notion very similar to FJ) in general and true KKT points {\color{\blk} under MFCQ or other related conditions (for example,~\cite{cartis2011evaluation,cartis2013evaluation,cartis2014complexity} assumed bounded penalty parameters in their problem settings, and the algorithm proposed in~\cite{birgin2016evaluation} obtained true KKT points whenever terminating at a certain algorithmic step)}. Further discussion of these related works in smooth optimization and other algorithmic approaches (beyond inexact proximal methods) are deferred to Section~\ref{relatedwork}.}

\subsection{Contribution}
We show that an inexact proximal method can solve a wide range of nonsmooth, {\color{black} weakly convex} constrained optimization problems, producing an approximate stationary point using at most $O(1/\epsilon^4)$ subgradient evaluations{\color{black}, matching the known rate for unconstrained optimization}. In particular, our proposed method uses a switching subgradient method approximately solving~\eqref{conssubproblem} to produce each subsequent $x_{k+1}$. Our analysis shows the following three generally desirable properties missing from prior works~\cite{ma2020quadratically,boob2022stochastic}:

\paragraph{Always Feasible Iterates}
By appropriately selecting the algorithmic parameters, we can ensure feasibility $ g_i(x_{k+1}) \leq 0$ at each iteration {\color{black} while maintaining a $O(1/\epsilon^{4})$ rate}. {\color{black} The method of~\cite{boob2022stochastic} can also maintain feasibility but at the cost of squaring the accuracy of subproblem solves needed, which would yield a strictly worse rate of convergence of $O(1/\epsilon^{6})$. The future development of more practical methods maintaining feasible iterates may be useful, for example, in settings of planning or control where feasibility corresponds to physical limitations or safety concerns~\cite{allgower2012nonlinear,yu2017design}. }

\paragraph{Stationarity with or without Constraint Qualification}
Ensuring constraint qualification over nonconvex constraints is nontrivial despite being consistently assumed by prior {\color{black} nonsmooth nonconvex} works. This is illustrated for a common sparse regularizer in Section~\ref{subsec:vignette} and numerical explored in Section~\ref{numerical}. In Theorems~\ref{theorem_mainFJ} and~\ref{theorem_mainKKT}, we show that at most $O(1/\epsilon^4)$ subgradient evaluations are required to produce an approximate KKT or FJ point, with or without constraint qualification, respectively.

\paragraph{Convergence Rates without Compactness}
Our guarantees apply without needing to assume compactness of the domain $X$, which prior works relied on. Hence, our theory applies more widely and, even in compact settings, may offer improvements as quantities like the diameter of $X$ are often replaced by smaller quantities dependent on the initialization. This is done by extending the analysis of the switching subgradient method to handle non-Lipschitz objective and constraint functions like those occurring in~\eqref{conssubproblem}. This analysis and resulting subproblem convergence guarantee appear to be new and may be of independent interest.

\subsection{Related Work}
\label{relatedwork}

{\color{black} \paragraph{Fritz-John and KKT Points in Smooth Optimization}
As mentioned in the introduction, much work has already been done in the context of smooth constrained optimization seeking FJ points in general and KKT points when constraint qualifications. A series of papers by Cartis, Gould, and Toint~\cite{cartis2011evaluation,cartis2013evaluation,cartis2014complexity} derived convergence guarantees for first and second-order methods for smooth nonconvex constrained problems. The penalized gradient method of~\cite{cartis2011evaluation} achieved a convergence rate of $O(1/\epsilon^2)$ towards KKT stationarity if constraint qualification holds {\color{\blk} (or more generally, whenever penalty parameters are bounded)} and in general, a slower $O(1/\epsilon^5)$ rate towards ``unscaled KKT points'', points whose violation of the KKT conditions is small in proportional to the size of the Lagrange multipliers. Subsequently, in~\cite{cartis2014complexity}, the authors showed a convergence rate of $O(1/\epsilon^2)$ is possible in both cases. This matches the rate for unconstrained smooth minimization, indicating no further progress is possible. In this work, we extend this message to nonsmooth settings, showing FJ and KKT convergence at a $O(1/\epsilon^4)$ rate matching the known rate for unconstrained {\color{black} weakly convex} nonsmooth minimization~\cite{davis2020nonsmooth}. Further discussions on scaled KKT points and needed bounds on Lagrange multipliers are given by~\cite{cartis2020strong}.

As wider advances in algorithm design, Hinder and Ye~\cite{hinder2018worst} showed that a (slightly modified) Fritz-John stationarity can be reached by an interior point method despite nonconvex constraints. See Bian et al.~\cite{bian2015ipm} as earlier work on nonconvex interior point method guarantees over simpler box constraints reaching scaled KKT points. Beyond just first and second-order methods, Birgin et al.~\cite{birgin2016evaluation} presented two-phase $p$th-order methods with guarantees towards producing KKT points. Their algorithms either terminate at scaled (or unscaled) KKT points or stop at an infeasible critical point. Beyond the first-order optimality conditions captured by FJ and KKT, {\color{\blk} second and higher-order} stationarity conditions have more recently been considered in{\color{\blk}~\cite{cartis2013evaluation,cartis2018second,cartis2020concise,nouiehed2018convergence,peng2020computation}}.}

\paragraph{Inexact Proximal Methods}
Using inexact proximal-point methods to solve nonsmooth {\color{black} weakly convex} problems is not new to this work. Double-loop algorithms that use several inner steps to inexactly solve a convex proximal subproblem in each outer iteration have been designed and analyzed widely. For example, the algorithm proposed in \cite{hare2009computing} approximating {\color{black} weakly convex} proximal points contributed to such an idea, and \cite{hare2010redistributed} presented a proximal variant of bundle methods solving such problems. More recently, \cite{davis2019proximally} developed this idea to give a $O(1/\epsilon^{4})$ convergence rate for unconstrained stochastic nonsmooth {\color{black} weakly convex} problems.

\paragraph{Special Case of (Strongly) Convex Constraints} 
A range of methods from the literature can be applied to inexactly solve the nonsmooth but strongly convex constrained subproblems constructed, which here arise as the subproblems~\eqref{conssubproblem}. A level-set method for structured convex constrained problems was introduced in \cite{aravkin2019level}, which was generalized and improved by \cite{lin2018level} to maintain feasibility. Alternative (augmented) Lagrangian approaches could be applied if near feasibility is sufficient.
Here, we take the approach of solving such problems via switching subgradient methods, which have been analyzed in \cite{articlepolyak} and extended in{\color{\blk}~\cite{ma2020quadratically,bayandina2018mirror,lan2020algorithms}}.

\paragraph{Comparison with Ma, Lin, and Yang~\cite{ma2020quadratically}}

We consider a very similar inexact proximal point method with switching subgradient method being the oracle for the subproblems as Ma et al. \cite{ma2020quadratically}, in which they also find nearly optimal and nearly feasible solutions for the subproblems. Their work also analyzed the convergence of a stochastic subgradient algorithm. However, in the deterministic setting, they only guarantee nearly feasible and approximate stationary solutions for the original optimization problem, while our method ensures actual feasibility. To attain KKT stationarity, they introduced a uniform Slater's condition as their stronger type of constraint qualification, which is stronger than our considered $\sigma$-strong MFCQ condition. Moreover, their upper bound on the optimal dual variables and convergence rates depend on the diameter of $X$, while we do not need such a requirement. Up to these constants, Ma et al.~proved a $O(1/\epsilon^4)$ rate of convergence towards KKT guarantees under MFCQ, which we match (in addition to our new FJ guarantees).

\paragraph{Comparison with Boob, Deng, and Lan~\cite{boob2022stochastic}}
As another closely related work, Boob et al.~\cite{boob2022stochastic} showed that the inexact proximal point method searching for nearly optimal solutions for the subproblems can ensure an approximate stationary solution for the main problem is found in a $O(1/\epsilon^4)$ convergence rate (assuming $X$ is compact). {\color{black} To additionally ensure feasibility, their framework requires a higher level of precision in subproblem solutions: namely, approximate solutions to~\eqref{conssubproblem} within distance $O(\epsilon^2)$ of the true solution are needed to maintain feasibility (see~\cite[(3.13)]{boob2022stochastic}), leading to an overall $O(1/\epsilon^6)$ rate. In contrast, our proposed method maintains feasibility while only requiring solutions within distance $O(\epsilon)$ of the true subproblem solution, giving an overall $O(1/\epsilon^4)$ rate.} Boob et al.~consider problem settings ranging from {\color{black} weakly convex} to strongly convex constrained problems and consider various MFCQ, strong MFCQ, and strong feasibility conditions as constraint qualifications. Their strong feasibility condition is stronger than our considered $\sigma$-strong MFCQ condition. Under their MFCQ and strong MFCQ conditions, an additional assumption is needed to ensure the existence of a stationary solution that the iterated points converge to and boundedness of the optimal dual variables. {\color{black} From strong MFCQ, we prove a constant upper bound on this directly.}

{\color{black}
\paragraph{Alternative Approaches given Nonconvex Constraints} Finally, we note several historically successful alternatives to (inexact) proximal methods: There is a long history of applying sequential quadratic programming techniques to smooth constrained optimization~\cite{Nocedal2006}. Of particular note,~\cite{curtis2012sequential} proposed and analyzed the limit points of such a method for locally Lipschitz nonsmooth nonconvex constrained problems. This was done without assuming weak convexity, making their method much more general than ours, but no explicit convergence rates were given. On the other hand, assuming smoothness, the cubic regularization approach of~\cite{Cartis2015}, penalized methods of~\cite{Facchinei2021,Wang2017} and augmented Lagrangian algorithm of~\cite{grapiglia2021complexity} can be applied with provably convergence guarantees, all faster than our $O(1/\epsilon^4)$ by virtue of their additional smoothness assumptions. If the constraints are star convex with respect to a known point (for example, the SCAD constraints previously mentioned), the radial methods of~\cite{radial1,radial2} could apply with the same $O(1/\epsilon^4)$ convergence rate we find while also maintaining feasible iterates.}

\subsection{Vignette: Failure of MFCQ Assumptions for Sparse Regularized Problems}
\label{subsec:vignette}
\begin{figure}[t]
    \centering
    \begin{subfigure}{0.25\textwidth}
    {
        \centering
        \includegraphics[width=\textwidth]{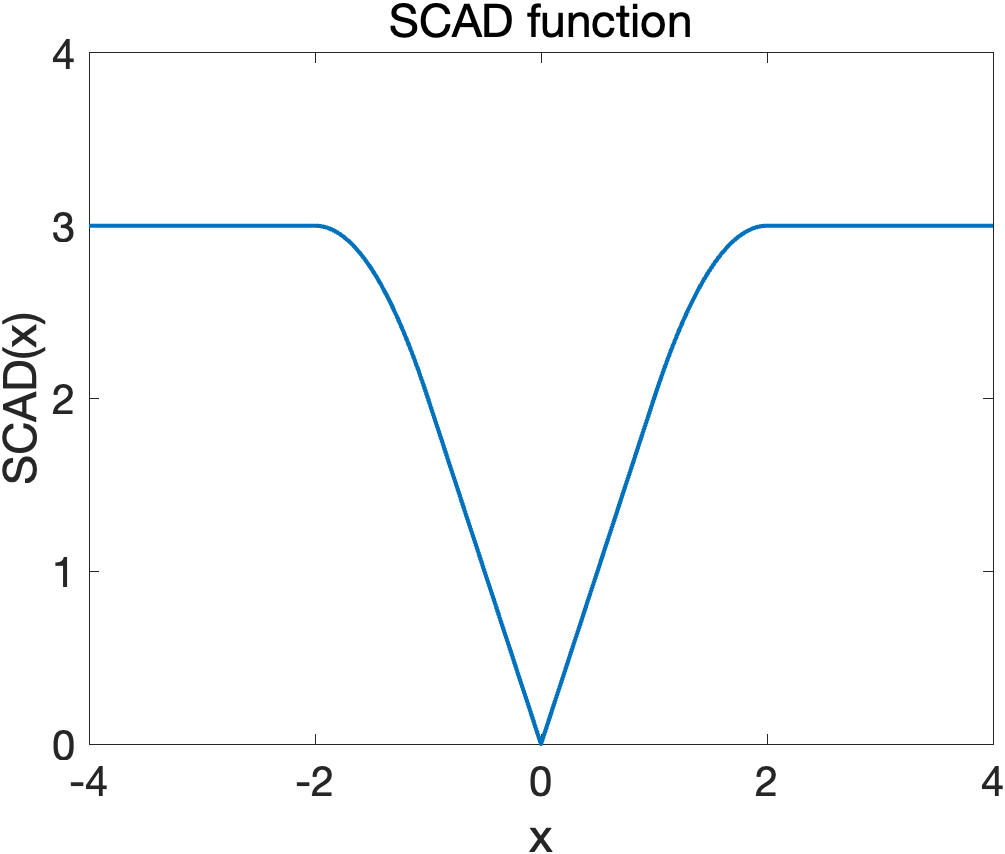}
        \caption{$1$D SCAD function}
        \label{1dSCAD}
    }
    \end{subfigure}
    \hfill
    \begin{subfigure}{0.70\textwidth}
    {
        \centering
        \includegraphics[width=\textwidth]{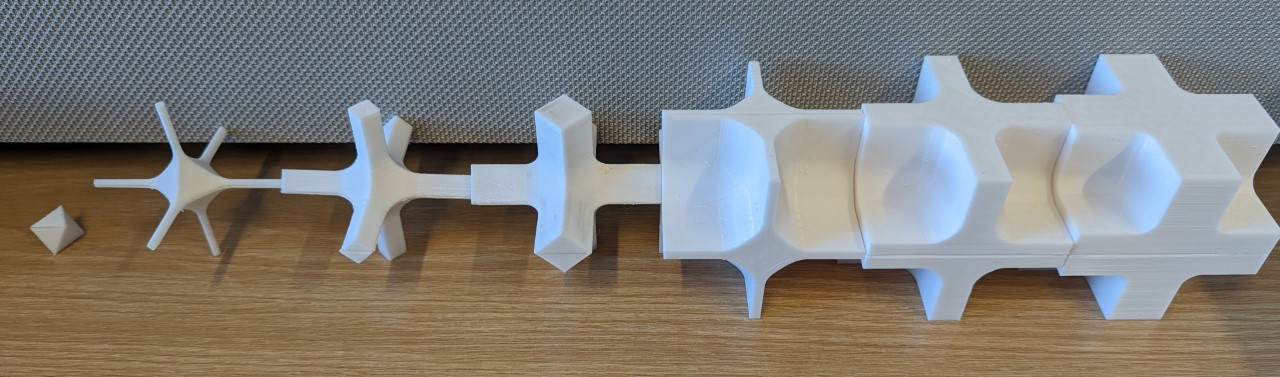}
        \caption{Seven SCAD level sets with $p\in\{2.5,3.5,4.5,5.5,6.5,7.5,8.5\}$. Note the set changes suddenly at $p\in\{3,6,9\}$, where MFCQ fails.}
        \label{3dSCAD}
    }
    \end{subfigure}
    \caption{The SCAD function $s$ and feasible regions in 3D given by $\sum_i s(x_i)\leq p$.}
    \label{scadplot}
\end{figure}

Nonsmooth nonconvex regularization has recently gained popularity due to its ability to facilitate stronger statistical guarantees on minimizers~\cite{wen2017efficient,wen2018survey,zhang2020meta,pieper2022nonconvex}.
One of the simplest regularizers is the Smoothly Clipped Absolute Deviation (SCAD) function~\cite{antoniadis1997wavelets,fan2001variable}, which sums up piecewise quadratic clipped absolute deviations in each coordinate
\begin{equation}
\begin{aligned}
    s(x_i)=
    \begin{cases}
        2|x_i| &{0 \leq |x_i| \leq 1}, \\
        -x_i^2+4|x_i|-1 &{1 < |x_i| \leq 2}, \\
        3 &{|x_i| > 2}.
    \end{cases}
    \label{scad}
\end{aligned}
\end{equation}
Near the origin, this behaves like a one-norm. As larger points are considered, it smoothly flattens out to avoid overly penalizing large entries. Figure~\ref{1dSCAD} shows the one-dimensional SCAD function. Note the constraint $g(x) := \sum_i s(x_i) - p \leq 0$ ensures that at most $\lfloor p/3 \rfloor$ entries of $x$ have a magnitude larger than two. Figure~\ref{3dSCAD} shows the feasible regions given by the three-dimensional SCAD constraints in $[-5,5]^3$. 

Optimization over these level sets will often yield sparse solutions. Since SCAD constraints are piecewise quadratic, we can often approximately solve the convex subproblem~\eqref{conssubproblem}. Despite this, two problems (one mild and one severe) prevent applying the convergence theory of prior works.

First, prior works do not apply as the set $\{x \mid g(x)\leq 0\}$ is not compact for any $p\geq 3$. If a bound on the size of a solution is known, then one could add a ball constraint $X = \{x \mid \|x\|\leq D\}$ to ensure compactness. Our theory applies without such a modification.

\begin{figure}[t]
    \centering
    \begin{subfigure}{0.69\textwidth}
    {
        \centering
        \includegraphics[width=\textwidth]{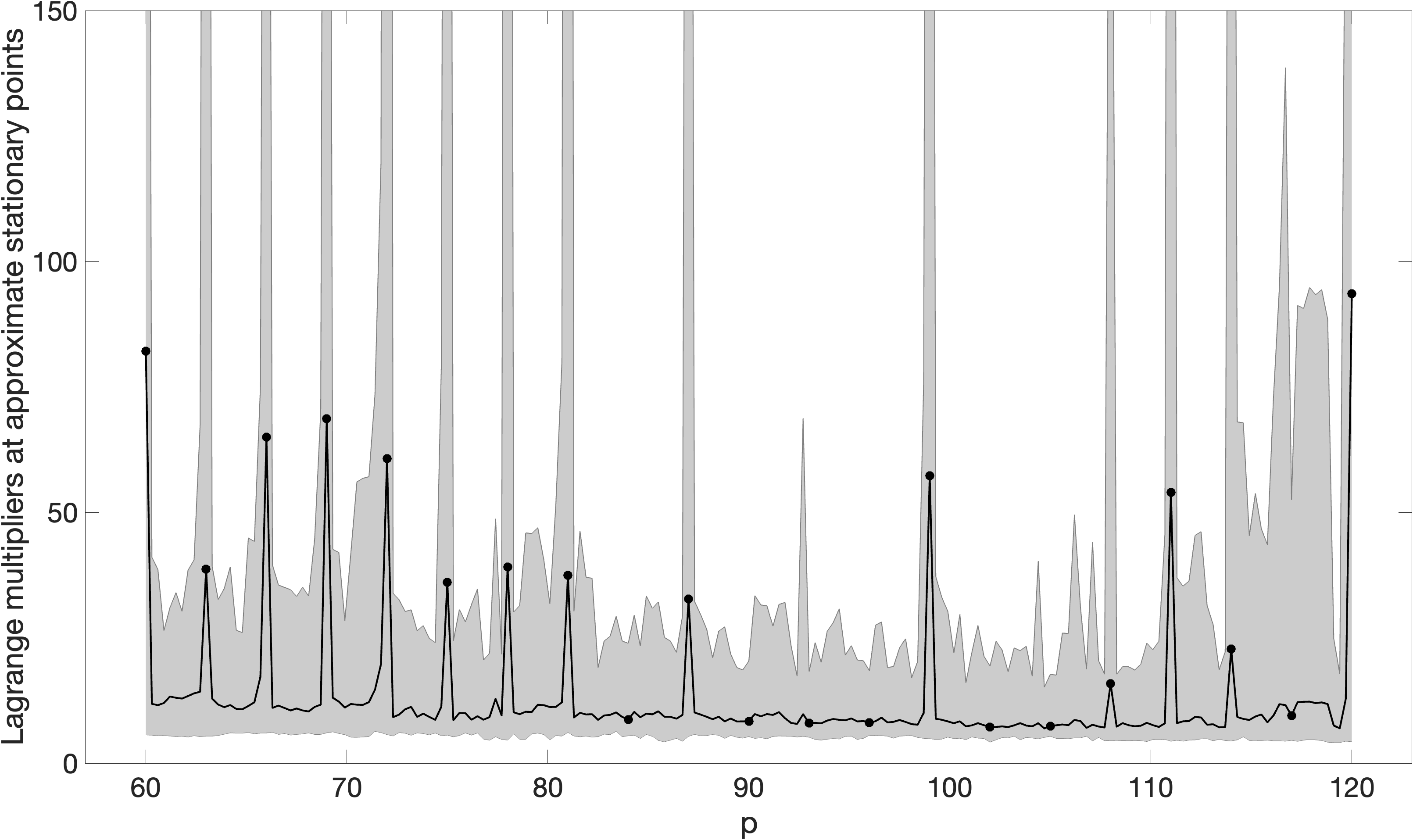}
        \caption{Small values of $p$}
    }
    \end{subfigure}
    \hfill
    \begin{subfigure}{0.251\textwidth}
    {
        \centering
        \includegraphics[width=\textwidth]{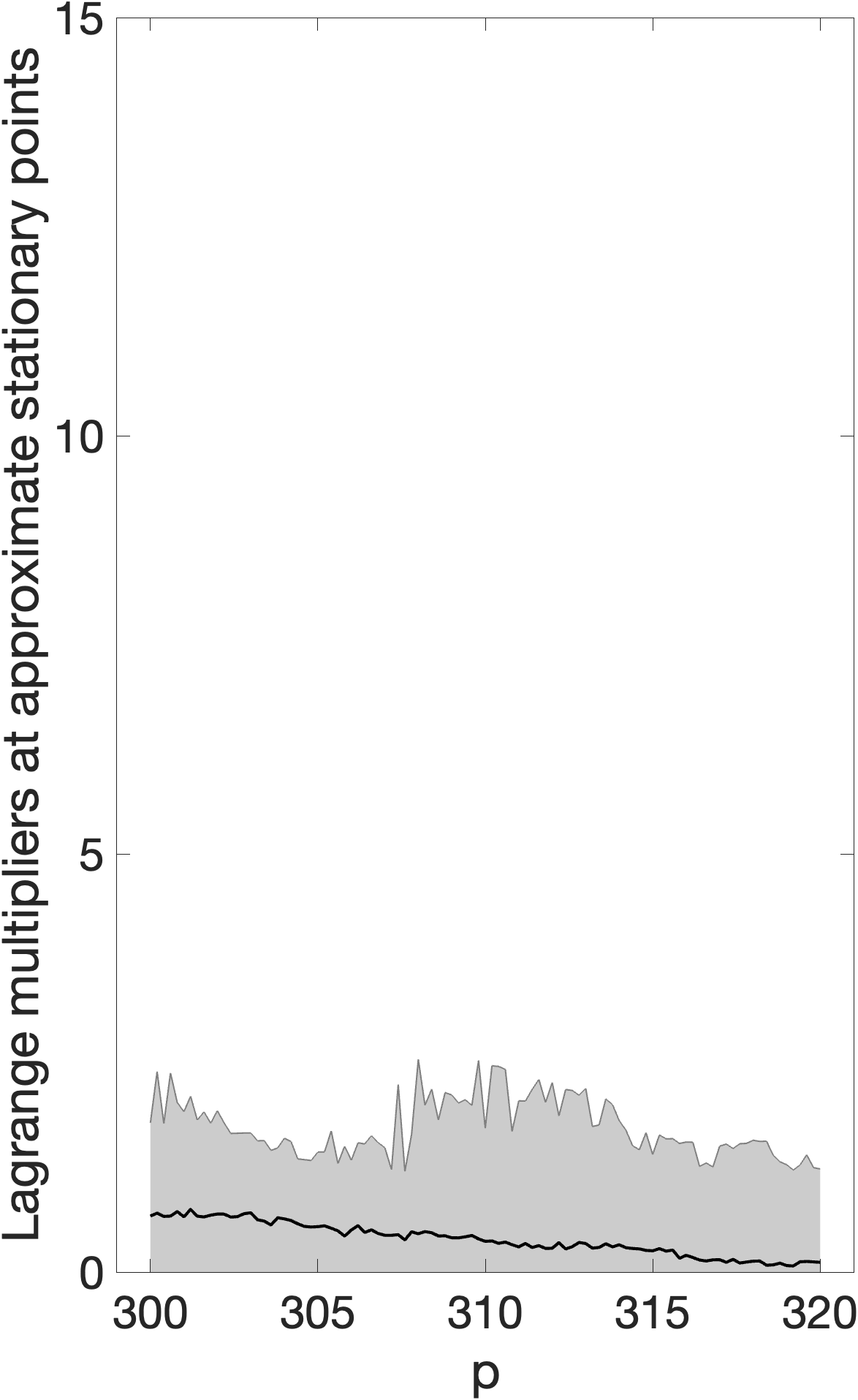}
        \caption{Large values of $p$}
    }
    \end{subfigure}
    \caption{Lagrange multipliers computed at approximate stationary points reached by iterating~\eqref{conssubproblem} on $30$ randomly generated SPR problems (see Section~\ref{numerical} for the exact construction). As $p$ varies from $60$ to $120$, the black line shows the average approximate multipliers reached and the gray region shows the range between maximum and minimum values seen. Black dots are placed at each multiple of three, where MFCQ fails to hold.}
    \label{example0}
\end{figure}

More subtly, prior works do not apply here as SCAD constraints often fail to have constraint qualification hold as $p$ varies. As a result, none of the prior works' theories provide any form of convergence guarantee. To illustrate this, we numerically consider the problem of Sparse Phase Retrieval problems (SPR), see \eqref{SPR}, which minimizes a piecewise quadratic objective over the piecewise quadratic constraint set for SCAD constraints. Figure~\ref{example0} shows the estimated Lagrange multipliers at limit points converged to by an inexact proximal point method. When $p$ is near a multiple of three, the limit point reached by iteratively applying~\eqref{conssubproblem} may fail to satisfy MFCQ, seen as its associated Lagrange multiplier blowing up, preventing KKT attainment. For large values of $p$, we see the multipliers tending to zero, corresponding to unconstrained stationarity.

Despite these failures of MFCQ, our theory still guarantees that the iteration will find an approximate FJ point. Note that Figure~\ref{example0} is based on averaging $30$ independent replicates. We only observe approximately {\color{\blk} $2\% \sim 10\%$} of replicates when $p$ is a multiple of three have their Lagrange multipliers diverge. So MFCQ is often violated but not everywhere. These sporadic failures in relatively simple nonsmooth nonconvex ({\color{black} weakly convex}) settings are one of this work's original motivations, leading us to develop theory capable of describing convergence when MFCQ fails while retaining (and improving) the convergence theory when MFCQ holds.

\section{Preliminaries}
\label{preliminary}

Throughout the paper, we use the following notations. Let $\|\cdot\|$ denote the $l_2$-norm. The distance from a point $x$ to a set $S$ is denoted as $\dist(x, S)=\min_{s \in S}\|x-s\|$, and the convex hull of any set $S$ is denoted as $\mathrm{co}\{S\}$. For any convex function $h:X \to \mathbb{R} \cup \{+\infty\}$, its set of subgradients at $x \in X$ is defined as:

\begin{equation}
\begin{aligned}
    \partial h(x)=\left\{\zeta \in \mathbb{R}^d|h(x') \geq h(x)+\zeta^T(x'-x), \quad \forall x' \in X\right\}.
\end{aligned}
\end{equation}

More generally, for any potentially nonconvex function $h:X \to \mathbb{R} \cup \{+\infty\}$, its set of Clarke subgradients at $x$ is defined as:

\begin{equation}
\begin{aligned}
    \partial h(x)=\mathrm{co}\left\{\lim_{i \to \infty}\nabla h(x_i)|x_i \to x \text{ and } h(x) \text{ is differentiable at any } x_i \in X\right\}.
\end{aligned}
\end{equation}
A function $h(x)$ is $\mu$-strongly convex on $X$ if $h-\frac{\mu}{2}\|\cdot\|^2$ is convex. This is equivalent to having:
\begin{equation}
\begin{aligned}
    h(x') \geq h(x)+\zeta^T(x'-x)+\frac{\mu}{2}\|x'-x\|^2, \qquad \forall x,x' \in X, \forall \zeta \in \partial h(x).
    \label{stronglyconvex}
\end{aligned}
\end{equation}
A function $h(x)$ is $\rho$-weakly convex on $X$ if $h+\frac{\rho}{2}\|\cdot\|^2$ is convex. This is equivalent to having:
\begin{equation}
\begin{aligned}
    h(x') \geq h(x)+\zeta^T(x'-x)-\frac{\rho}{2}\|x'-x\|^2, \qquad \forall x,x' \in X, \forall \zeta \in \partial h(x).
    \label{weaklyconvex}
\end{aligned}
\end{equation}

We consider two different notions describing approximate stationarity for our nonsmooth nonconvex constrained problem of interest~\eqref{mainproblemoriginal}, weakening the FJ conditions and KKT conditions shown in \eqref{FJoriginal} and \eqref{KKToriginal} respectively.

\begin{definition}\label{definitionFJ1}
    A point x is an $\epsilon$-FJ point for problem \eqref{mainproblemoriginal} if $g_i(x) \leq 0 \ \forall i=1,...,m$, and there exists $\zeta_f \in \partial f(x)$, $\zeta_{gi} \in \partial g_i(x)$ and $\gamma_0 \geq 0$, $\gamma=(\gamma_1,...,\gamma_m)^T \geq 0$, $\gamma_0+\sum_{i=1}^m\gamma_i=1$ such that:
    \begin{align}
        \dist\left(\gamma_0\zeta_f+\sum_{i=1}^m\gamma_i\zeta_g, -N_X(x)\right) &\leq \epsilon,\label{FJ1_1} \\
        |\gamma_i g_i(x)| &\leq \epsilon^2 \qquad \forall i=1,...,m.\label{FJ1_2}
    \end{align}
\end{definition}

\begin{definition}\label{definitionKKT1}
    A point x is an $\epsilon$-KKT point for problem \eqref{mainproblemoriginal} if $g_i(x) \leq 0 \ \forall i=1,...,m$, and there exists $\zeta_f \in \partial f(x)$, $\zeta_{gi} \in \partial g_i(x)$ and $\lambda=(\lambda_1,...,\lambda_m)^T \geq 0$ such that:
    \begin{align}
        \dist\left(\zeta_f+\sum_{i=1}^m\lambda_i\zeta_{gi}, -N_X(x)\right) &\leq \epsilon, \label{KKT1_1}\\
        |\lambda_i g_i(x)| &\leq \epsilon^2 \qquad \forall i=1,...,m.\label{KKT1_2}
    \end{align}
\end{definition}

Let $\hat{x}_{k+1}$ denote the optimal solution for the subproblem \eqref{conssubproblem}. The considered inexact proximal point approach will produce iterates $x_{k+1}$ near each $\hat{x}_{k+1}$. As we will see, the sequence $\hat{x}_{k}$ converges towards an approximate stationary point for the main problem \eqref{mainproblemoriginal}. So we can only ensure our iterates $x_k$ are near an approximately stationary point. The following definitions describe points in the proximity of an approximately stationary point.

\begin{definition}
    A {\color{black} feasible} point x is an $(\epsilon, \eta)$-FJ point for problem \eqref{mainproblemoriginal} if there exists an $\epsilon$-FJ point $x'$ for problem \eqref{mainproblemoriginal} with $\|x-x'\| \leq \eta$.
    \label{FJpoint2}
\end{definition}

\begin{definition}
    A {\color{black} feasible} point x is an $(\epsilon, \eta)$-KKT point for problem \eqref{mainproblemoriginal} if there exists an $\epsilon$-KKT point $x'$ for problem \eqref{mainproblemoriginal} with $\|x-x'\| \leq \eta$.
    \label{KKTpoint2}
\end{definition}

The accuracy of KKT stationarity guarantees we derive will depend on the sizes of the associated Lagrange multipliers. To give a constant upper bound on these optimal Lagrange multipliers in \eqref{KKToriginal} for our subproblems (see problem \eqref{subproblem} below), we assume a stronger type of constraint qualification defined below. Let $A(x)=\{i \mid g_i(x)=0, i=1,...,m\}$. We say $\sigma$-strong MFCQ condition holds at $x$ if there exists a constant $\sigma>0$, such that:

\begin{equation}
\begin{aligned}
    \exists v \in -N_X^*(x) \quad and \quad \|v\|=1 \qquad s.t. \quad \zeta_{gi}^Tv \leq -\sigma \quad \forall i \in A(x), \forall \zeta_{gi} \in \partial g_i(x).
    \label{strongMFCQ}
\end{aligned}
\end{equation}

We say the $\sigma$-strong MFCQ condition holds for problem~\eqref{mainproblemoriginal} when $\sigma$-strong MFCQ condition is satisfied at any $x \in X$. When the $\sigma$-strong MFCQ condition is satisfied for all the subproblems~\eqref{subproblem}, our Lemma~\ref{lemmadualvar} shows boundedness of Lagrange multipliers in \eqref{KKToriginal} for our subproblems. This boundness is critical to improve our FJ convergence guarantees to convergence towards KKT stationarity.

Without loss of generality, we simplify the $m$ nonsmooth, {\color{black} weakly convex} constraints of~\eqref{mainproblemoriginal} into a single constraint as follows:
\begin{equation}
\begin{cases}
    \min_{x \in X} \quad &f(x) \\
    \mathrm{s.t.} \quad &g(x):=\max_{i=1,...,m}g_i(x) \leq 0.
    \label{mainproblem}
\end{cases}
\end{equation}
Note if each $g_i$ is $\rho$-weakly convex, then $g$ is $\rho$-weakly convex. Note that in this reformulation, there is only a single constraint and hence only a single Lagrange multiplier. Since subgradients of a finite maximum of $m$ elements are convex combinations of subgradients of the component functions, the original vector of multipliers can always be recovered.

We follow the same construction as~\eqref{conssubproblem} to build our proximal subproblems for~\eqref{mainproblem}, given by

\begin{equation}
\begin{cases}
    \min_{x \in X} \quad &F_k(x):=f(x)+\frac{\hat{\rho}}{2}\|x-x_k\|^2 \\
    \mathrm{s.t.} \quad &G_k(x):=g(x)+\frac{\hat{\rho}}{2}\|x-x_k\|^2 \leq 0.
    \label{subproblem}
\end{cases}
\end{equation}
By selecting $\hat{\rho} > \rho$, both the objective function $F_k(x)$ and the constraint $G_k(x)$ are $(\hat{\rho}-\rho)$-strongly convex. Throughout, we require $\hat{\rho}>\max\{\rho,1\}$. In its outer loop, our inexact proximal point method will set $x_{k+1}$ as a nearly optimal and feasible solution of~\eqref{subproblem}.

We make the following four assumptions about \eqref{mainproblem} throughout this paper.

\begin{assumption}
    $f(x)$ and $g(x)$ are continuous and $\rho$-weakly convex on $X$.
    \label{assumption1}
\end{assumption}
\begin{assumption}
    $f_{lb}=\inf_{x \in X}f(x)>-\infty$, $g_{lb}=\inf_{x \in X}g(x)>-\infty$.
    \label{assumption2}
\end{assumption}
\begin{assumption}
    For any $x \in X$, we can compute $\zeta_f \in \partial f(x)$, $\zeta_g \in \partial g(x)$ with $\|\zeta_f\|, \|\zeta_g\| \leq M$.
    \label{assumption3}
\end{assumption}
\begin{assumption}
    A feasible point $x_0$ is known (i.e. $x_0 \in X$ and $g(x_0) \leq 0$).
    \label{assumption4}
\end{assumption}

These assumptions suffice for our convergence theory to FJ points. Under the following additional assumption, our convergence results improve to ensure approximate KKT stationarity.

\begin{assumption}
    $\sigma$-strong MFCQ condition is satisfied for any subproblem \eqref{subproblem}.
    \label{assumption5}
\end{assumption}

Let $\hat{x}_{k+1}$ denote the optimal solution for the subproblem \eqref{subproblem}. In the following lemma, we will show that when $\|\hat{x}_{k+1}-x_k\|$ is small enough, the first conditions for either FJ or KKT stationarity \eqref{FJ1_1}/\eqref{KKT1_1} hold for the original nonsmooth {\color{black} weakly convex} problem~\eqref{mainproblem}. Further utilizing the selection $\hat{\rho}>\max\{\rho,1\}$, we conclude the second conditions \eqref{FJ1_2}/\eqref{KKT1_2} must be satisfied when \eqref{FJ1_1}/\eqref{KKT1_1} are. Hence our convergence theory follows along the following reasoning: once $\hat{x}_{k+1}$ is an approximate stationary point for the main problem \eqref{mainproblem}, $x_k$ must lie in a neighborhood of $\hat{x}_{k+1}$. Then depending on whether Assumption E holds, this gives an approximate FJ or KKT stationary solution near $x_k$.

\begin{lemma}
    When Assumptions \ref{assumption1}-\ref{assumption4} hold and $\hat{\rho}>\max\{\rho,1\}$, if $\|\hat{x}_{k+1}-x_k\| \leq \frac{\epsilon}{\hat{\rho}}$ then $x_k$ is an $(\epsilon,\epsilon)$-FJ point. If additionally, Assumptions \ref{assumption5} holds, then a dual optimal $\lambda_k$ for~\eqref{subproblem} exists and if $\|\hat{x}_{k+1}-x_k\| \leq \frac{\epsilon}{\hat{\rho}(1+\lambda_k)}$, then $x_k$ is an $(\epsilon,\epsilon)$-KKT point.
    \label{lemmastationarity}
\end{lemma}

Note in the second case (under Assumption~\ref{assumption5}), the size  of the Lagrange multiplier plays a role. As $\lambda_k$ grows larger, the stationarity $\|\hat{x}_{k+1}-x_k\|$ needs to be smaller to ensure the same level of KKT attainment. For notational ease, let $D=\sqrt{(-8g_{lb})/(\hat{\rho}-\rho)}$ denote the upper bound on the diameter of every subproblem constraint set $\{x \mid G_k(x) \leq 0\}$ due to the $(\hat{\rho}-\rho)$-strong convexity of $G_k(x)$\footnote{{\color{black} For completeness, we provide a short verification of this: Let $x^*$ be the minimizer of $G_k$. Then any $x_i$ with $G_k(x_i) \leq 0$ have by $(\hat{\rho}-\rho)$-strong convexity $0 \geq G_k(x_i) \geq G_k(x^*)+\frac{\hat{\rho}-\rho}{2}\|x_i-x^*\|^2 \geq g_{lb}+\frac{\hat{\rho}-\rho}{2}\|x_i-x^*\|^2$. Therefore $\|x_1-x_2\| \leq \|x_1-x^*\|+\|x_2-x^*\| \leq \sqrt{(-8g_{lb})/(\hat{\rho}-\rho)}$.}}. In particular, this upper bounds the distance from the current iterate $x_k$ to $\hat x_{k+1}$. Using this, in Lemma~\ref{lemmadualvar}, we show that $\sigma$-strong MFCQ (Assumption~\ref{assumption5}) ensures a uniform upper bound for the optimal subproblem dual variables of $B=\frac{M+\hat{\rho}D}{\sigma}$.
    \section{Algorithms}

This section first describes the switching subgradient method and, second, our use of it as an oracle for solving the main problem \eqref{mainproblem} in our inexact proximal point method. All proofs are deferred to Section \ref{analysis}. {\color{black} The switching subgradient method is not fundamentally needed. Any subroutine would suffice if it is capable of minimizing a Lipschitz function $f$ plus a quadratic subject to several Lipschitz functions $g_i$ plus a quadratic being nonpositive. If one can project onto the subproblem's constraint set efficiently, a projected subgradient method could be applied (Remark~\ref{rem:proj-subgrad} discusses such a procedure for the special case of SCAD constraints). Any proximal subgradient-type method able to handle functional constraints would also be a reasonable alternative, but we are unaware of any method that has sufficient theoretical support to meet our needs.}

\subsection{The Classic Switching Subgradient Method (without Lipschitz Continuity)}

\label{switchingsubgradient}
We introduce the classic switching subgradient method (see \cite{articlepolyak}) for solving problems of the form
\begin{equation}
\begin{cases}
    \min_{z \in Z} \quad &F(z) \\
    \mathrm{s.t.} \quad &G(z) \leq 0.
    \label{generalsubproblem}
\end{cases}
\end{equation}
We assume the domain $Z$ is a closed convex set, and $F(z)$ and $G(z)$ are $\mu$-strongly convex functions on $Z$. Let $z^*$ be the unique optimal solution to this problem.
We define nearly optimal and nearly feasible solutions for this problem as follows.
\begin{definition}
    A point $z$ is a $(\delta,\tau)$-optimal solution for problem \eqref{generalsubproblem} if $F(z)-F(z^*) \leq \delta$ and $G(z) \leq \tau$, where $z^*$ is the optimal solution.
\end{definition}

Here, we analyze the switching subgradient method (Algorithm \ref{algorithm_sub}) to solve problem \eqref{generalsubproblem}, finding a $(\tau,\tau)$-optimal solution for it. When the current iterate is not nearly feasible with tolerance $\tau$, we compute the subgradient based on the constraint function and make an update seeking feasibility; otherwise, we compute the subgradient of the objective function to make an update seeking optimality.

\begin{algorithm}
\caption{The Switching Subgradient Method $SSM(\tau,T,z_0,\{\alpha_t\})$}
\label{algorithm_sub}
\begin{algorithmic}
\REQUIRE{$\tau>0$, $T>0$, $z_0\in Z$, $\{\alpha_t\}_{t=0}^{T-1}$}
\STATE{Set $I=\emptyset$, $J=\emptyset$}
\FOR{$t=0,1,...,T-1$}
\IF{$G(z_t) \leq \tau$}
\STATE{$z_{t+1}=\proj_Z(z_t-\alpha_t\zeta_{Ft})$, $\zeta_{Ft} \in \partial F(z_t)$, $I=I \cup \{t\}$}
\ELSE
\STATE{$z_{t+1}=\proj_Z(z_t-\alpha_t\zeta_{Gt})$, $\zeta_{Gt} \in \partial G(z_t)$, $J=J \cup \{t\}$}
\ENDIF
\ENDFOR
\ENSURE{$\Bar{z}_T = \frac{\sum_{t \in I}(t+1){\color{\blk} z_t}}{\sum_{t \in I}(t+1)}$}
\end{algorithmic}
\end{algorithm}

We give the convergence result for this method, generalizing \cite{lan2020algorithms,bayandina2018mirror,ma2020quadratically}.
These previous convergence analyses have assumed uniform Lipschitz continuity for both $F(z)$ and $G(z)$. However, such results are insufficient for analyzing its application to~\eqref{subproblem} since the added quadratic terms rule out global Lipschitz continuity. Instead, for our analysis here, we only need the following weaker, non-Lipschitz condition, previously considered for projected subgradient methods~\cite{grimmer2019convergence}: For any given target level of feasibility $\tau$, suppose there exist constants $L_0, L_1 \geq 0$ such that all nearly feasible $z_1 \in \{z \mid G(z)\leq \tau\}$ and infeasible $z_2\in \{z \mid G(z) >\tau\}$ have subgradients $\zeta_F \in \partial F(z_1),\ \zeta_G \in \partial G(z_2)$ bounded affinely by their current suboptimality/infeasibility
\begin{equation}
\begin{aligned}
    \|\zeta_F\|^2 &\leq L_0^2+L_1(F(z_1)-F(z^*)),\\
    \|\zeta_G\|^2 &\leq L_0^2+L_1(G(z_2)-G(z^*)).
    \label{quadraticbound}
\end{aligned}
\end{equation}
When $L_1=0$, this captures the standard case of $L_0$-Lipschitz $F(z)$ and $G(z)$. However, no function can possess Lipschitz continuity and strong convexity on an unbounded domain. When $L_1>0$, the non-Lipschitz condition~\eqref{quadraticbound} allows $F(z)$ and $G(z)$ to grow quadratically (hence this assumption is not at odds with strong convexity on unbounded domains).

\begin{theorem} \label{theorem_generalsub}
    Given $\alpha_t=\frac{2}{\mu(t+2)+\frac{L_1^2}{\mu(t+1)}}$, $\tau>0$, and {\color{\blk}$z_0\in Z$} with $G(z_0)\leq \tau$, Algorithm \ref{algorithm_sub}'s output $\Bar{z}_T$ is a $(\tau,\tau)$-optimal solution for problem \eqref{generalsubproblem} for all
    \begin{equation}
        T \geq \max\left\{\frac{8L_0^2}{\mu\tau},\sqrt{\frac{2L_1^2\|z_0-z^*\|^2}{\mu\tau}}\right\}.\nonumber
    \end{equation}
\end{theorem}
Minor modifications of our analysis would show that the switching subgradient method can attain a $(\tau,0)$-optimal solution at the rate of $O(1/\tau)$ for problem \eqref{generalsubproblem}, provided a strictly feasible Slater point exists (i.e., $G(z_0) < 0$).

In the proximal subproblem~\eqref{subproblem}, $F_k$ and $G_k$ are both $(\hat{\rho}-\rho)$-strongly convex functions. Consequently, they are not Lipschitz if the domain $X$ is unbounded. In the following lemma, however, we bound its subgradients via the non-Lipschitz condition~\eqref{quadraticbound}. Guarantees for the switching subgradient method applied to these proximal subproblems directly follow.
\begin{lemma}
    For any $x_k\in X$ with $g(x_k)\leq 0$, the non-Lipschitz condition~\eqref{quadraticbound} is satisfied by the proximal subproblem~\eqref{subproblem} with $L^2_0=9M^2-6\hat{\rho}g_{lb}$ {\color{\blk} and} $L_1=6\hat{\rho}$.
    \label{lemma_quadraticbound}
\end{lemma}

\begin{corollary}\label{corollary_sub}
    With $z_0=x_k$, $\alpha_t=\frac{2}{(\hat{\rho}-\rho)(t+2)+\frac{36\hat{\rho}^2}{(\hat{\rho}-\rho)(t+1)}}$ and $\tau>0$ in Algorithm \ref{algorithm_sub}, $\Bar{z}_T$ is a $(\tau,\tau)$-optimal solution for problem \eqref{subproblem} for all
    \begin{equation}
        T \geq \max\left\{\frac{24(3M^2-2\hat{\rho}g_{lb})}{(\hat{\rho}-\rho)\tau},\sqrt{\frac{72\hat{\rho}^2D^2}{(\hat{\rho}-\rho)\tau}}\right\}.\nonumber
    \end{equation}
\end{corollary}

In prior convergence analysis of the switching subgradient method, Lipschitz continuity is necessary for both the objective function $F_k(x)$ and the constraint function $G_k(x)$. Since these functions are strongly convex (and so grow quadratically), previous works required compactness of the domain $X$ to yield a uniform Lipschitz constant. In contrast, our Corollary \ref{corollary_sub} avoids assuming any compactness.

Several stochastic variants of Algorithm \ref{algorithm_sub} have been considered for solving stochastic generalizations of~\eqref{generalsubproblem}. An adaptive stochastic mirror descent method was introduced in~\cite{bayandina2018mirror}, which assumes exact functional values are computable for each constraint, but only stochastic approximations of the subgradients of the objective and constraints are available. With unbiased estimators of the subgradients, Algorithm \ref{algorithm_sub} can be applied to this kind of stochastic problem with convergence results in expectation without requiring the compactness of the domain or the stochastic subgradients to be bounded almost surely. A stochastic version of the non-Lipschitz condition~\eqref{quadraticbound} was considered by \cite{gower2019sgd} as a combination of the expected smoothness and finite gradient noise conditions around the optimal solution, which is needed to show convergence of the stochastic switching subgradient method. In \cite{lan2020algorithms}, they proposed a cooperative stochastic approximation method under stochastic estimations of the functional values of both the objective function and the constraint. Under this setting, they showed guarantees of finding nearly optimal solutions in expectation (although still requiring the compactness of the domain).

{\color{black}
\begin{remark} \label{rem:proj-subgrad}
    If one extends our computational model to assume projections onto the constraint set $\{z\in Z \mid G(z)\leq 0\}$ are tractable, a simpler projected subgradient method can be employed than the above switching subgradient method:
    $$ z_{t+1} = \mathrm{proj}_{\{z\in Z \mid G(z)\leq 0\}}(z_t - \alpha_t\zeta_{Ft}), \ \zeta_{Ft}\in \partial F(z_t) \ . $$
    The analysis of~\cite[Theorem 1.7]{grimmer2019convergence} shows such an iteration converges at the same rate as Theorem~\ref{theorem_generalsub}. Consequently, our theory applies equally when using a projected subgradient method to solve proximal subproblems, provided this stronger projection oracle is available. For example, in the motivating setting of SCAD constraints in Section~\ref{subsec:vignette} with $Z=\mathbb{R}^n$, this projection can be done efficiently, see Appendix~\ref{app:SCAD_proj} for the details of such a calculation. Although such specialized settings are not our primary focus, a projection-based subroutine there would yield actually sparse iterates.
\end{remark}
}

\subsection{Proximally Guided Switching Subgradient Method}

Our primary method of interest iteratively uses the switching subgradient method to inexactly produce proximal point steps, following the idea of~\eqref{conssubproblem}. This process of repeatedly approximately solving~\eqref{subproblem} is formalized in Algorithm~\ref{algorithm_main}.

\begin{algorithm}
\caption{The Proximally Guided Switching Subgradient Method}
\begin{algorithmic}
\REQUIRE{$\hat{\rho}>\max\{\rho,1\}$, $\tau>0$, $T_{inner}$, $x_0 \in X$ with $g(x_0)\leq 0$.}
\STATE Set $\alpha_t=\frac{2}{(\hat{\rho}-\rho)(t+2)+\frac{36\hat{\rho}^2}{(\hat{\rho}-\rho)(t+1)}}$
\FOR{$k=0,1,...,$}
\STATE{Set $x_{k+1}$ as the output of $SSM(\tau,T_{inner},x_k,\{\alpha_t\})$ applied to~\eqref{subproblem}}
\ENDFOR
\end{algorithmic}
\label{algorithm_main}
\end{algorithm}
Our primary result is that this simple scheme will produce Fritz-John points whenever the Assumptions \ref{assumption1}--\ref{assumption4} hold (amounting to standard bounds on continuity, nonconvexities, objective values, and the initialization). When constraint qualification (via Assumption~\ref{assumption5}) is additionally assumed, our theory improves to ensure a KKT point is found. To derive {\color{\blk} these} improved approximate KKT guarantees, we show that this additional assumption yields a uniform upper bound for the optimal dual variables (Lagrange multipliers) of the KKT conditions \eqref{KKToriginal} for each of the subproblems~\eqref{subproblem}. This is formalized in the following lemma.
\begin{lemma}
    Under Assumptions \ref{assumption1}--\ref{assumption5}, the optimal dual variables for problems \eqref{subproblem} are uniformly upper bounded by $B=\frac{M+\hat{\rho}D}{\sigma}$.
    \label{lemmadualvar}
\end{lemma}

To guarantee the identification of an $(\epsilon, \epsilon)$-FJ point or $(\epsilon, \epsilon)$-KKT point, our theory requires slightly different selections for the feasibility tolerance $\tau$ and how many iterations $T_{inner}$ of the inner switching subgradient method to utilize. Namely, in these two different settings respectively, we select
\begin{align}
    &\begin{cases}
        \tau_{FJ}=\frac{(\hat{\rho}-\rho)\epsilon^2}{4\hat{\rho}(2\hat{\rho}-\rho)}\\
        T_{FJ} = \max\left\{\frac{96\hat{\rho}(2\hat{\rho}-\rho)(3M^2-2\hat{\rho}g_{lb})}{(\hat{\rho}-\rho)^2\epsilon^2},\sqrt{\frac{288\hat{\rho}^3(2\hat{\rho}-\rho)D^2}{(\hat{\rho}-\rho)^2\epsilon^2}}\right\} \ ,
    \end{cases}
    \label{paraFJ}\\
    &\begin{cases}
        \tau_{KKT}=\frac{(\hat{\rho}-\rho)\epsilon^2}{4(1+B)^2\hat{\rho}(2\hat{\rho}-\rho)}\\
        T_{KKT} = \max\left\{\frac{96(1+B)^2\hat{\rho}(2\hat{\rho}-\rho)(3M^2-2\hat{\rho}g_{lb})}{(\hat{\rho}-\rho)^2\epsilon^2},\sqrt{\frac{288(1+B)^2\hat{\rho}^3(2\hat{\rho}-\rho)D^2}{(\hat{\rho}-\rho)^2\epsilon^2}}\right\}\ .
    \end{cases}
    \label{paraKKT}
\end{align}
These feasibility tolerances are chosen as they guarantee the feasibility of the iterates $x_k$ of Algorithm \ref{algorithm_main} until an appropriate FJ or KKT point is found. This is formalized in the following lemma.
\begin{lemma}
    Under Assumptions \ref{assumption1}--\ref{assumption4} with $\tau$ and $T_{inner}$ as in \eqref{paraFJ} (or under Assumptions \ref{assumption1}--\ref{assumption5} with $\tau$ and $T_{inner}$ as in \eqref{paraKKT}) {\color{black} and any $\hat\rho >\max\{\rho,1\}$,} Algorithm \ref{algorithm_main} has $g(x_k) \leq 0$ at every iteration $k$ before $\hat x_{k+1}$ is an $\epsilon$-FJ point  (or an $\epsilon$-KKT point).
    \label{lemmafea}
\end{lemma}

As a result, one need not worry about the proposed method becoming infeasible and converging to a stationary point outside the feasible region. The following pair of theorems then guarantee that at most {\color{\blk} $O(1/\epsilon^4)$} subgradient evaluations are needed for this feasible sequence of iterates to reach an approximate FJ or KKT point.
\begin{theorem}
    Under Assumptions \ref{assumption1}--\ref{assumption4} and any $\hat\rho >\max\{\rho,1\}$, Algorithm \ref{algorithm_main} with $\tau = \tau_{FJ}$ and $T_{inner}=T_{FJ}=:\max\left\{\frac{\Delta_1}{\epsilon^2},\frac{\Delta_2}{\epsilon}\right\}$ has $x_K$ be an $(\epsilon, \epsilon)$-FJ point for problem \eqref{mainproblem} for some
    \begin{equation}
        K \leq \frac{4\hat{\rho}^2(f(x_0)-f_{lb})}{(\hat{\rho}-\rho)\epsilon^2}=:\frac{\Delta_3}{\epsilon^2} \ . \nonumber
    \end{equation}
    Such an $x_K$ is found using at most $\frac{\Delta_3\max\left\{\Delta_1,\Delta_2\epsilon\right\}}{\epsilon^4}$
    total subgradient evaluations\footnote{{\color{black} The definitions of the constants $\Delta_1$, $\Delta_2$ and $\Delta_3$ can be found at the end of the proof.}}.
    \label{theorem_mainFJ}
\end{theorem}

\begin{theorem}
    Under Assumptions \ref{assumption1}--\ref{assumption5} and any $\hat\rho >\max\{\rho,1\}$, Algorithm \ref{algorithm_main} with $\tau = \tau_{KKT}$ and $T_{inner}=T_{KKT}=:\max\left\{\frac{\Lambda_1}{\epsilon^2},\frac{\Lambda_2}{\epsilon}\right\}$ has $x_K$ be an $(\epsilon, \epsilon)$-KKT point for problem \eqref{mainproblem} for some
    \begin{equation}
        K \leq \frac{4(1+B)\hat{\rho}^2(f(x_0)-f_{lb})}{(\hat{\rho}-\rho)\epsilon^2}=:\frac{\Lambda_3}{\epsilon^2} \ . \nonumber
    \end{equation}
    Such an $x_K$ is found using at most $\frac{\Lambda_3\max\left\{\Lambda_1,\Lambda_2\epsilon\right\}}{\epsilon^4}$
    total subgradient evaluations\footnote{{\color{black} The definitions of the constants $\Lambda_1$, $\Lambda_2$ and $\Lambda_3$ can be found at the end of the proof.}}.
    \label{theorem_mainKKT}
\end{theorem}

\section{Convergence Analysis}
\label{analysis}

\subsection{Proof of Theorem \ref{theorem_generalsub}}
Our convergence proof for the switching subgradient method presented here follows closely in the styles of~\cite{lacostejulien2012simpler,ma2020quadratically,grimmer2019convergence}.
Let $z^*$ be the optimal solution for~\eqref{generalsubproblem}, whose existence and uniqueness follow from strong convexity. When $t \in I$, we have
\begin{align}
    \|z_{t+1}-z^*\|^2 &\leq \|z_t-\alpha_t\zeta_{Ft}-z^*\|^2 \nonumber\\
    &=\|z_t-z^*\|^2-2\alpha_t\zeta_{Ft}^T(z_t-z^*)+\alpha_t^2\|\zeta_{Ft}\|^2 \nonumber\\
    &\leq \|z_t-z^*\|^2-2\alpha_t\zeta_{Ft}^T(z_t-z^*)+L_0^2\alpha_t^2+L_1\alpha_t^2(F(z_t)-F(z^*))\nonumber\\
    &\leq (1-\mu\alpha_t)\|z_t-z^*\|^2-(2\alpha_t-L_1\alpha_t^2)(F(z_t)-F(z^*))+L_0^2\alpha_t^2 \ . \nonumber
\end{align}
where the first inequality uses the nonexpansiveness of projections, the second uses the non-Lipschitz subgradient bound, and the third uses strong convexity.
Hence 
\begin{align}
    (2-L_1\alpha_t)(F(z_t)-F(z^*)) &\leq \left(\frac{1}{\alpha_t}-\mu\right)\|z_t-z^*\|^2-\frac{1}{\alpha_t}\|z_{t+1}-z^*\|^2+L_0^2\alpha_t\ .\nonumber
\end{align}
Since $\alpha_t=\frac{2}{\mu(t+2)+\frac{L_1^2}{\mu(t+1)}}$, the above coefficient on $F(z_t)-F(z^*)$ is at least one, i.e., 
\begin{equation}
    L_1\alpha_t=\frac{2L_1}{\mu(t+2)+\frac{L_1^2}{\mu(t+1)}} \leq \frac{2L_1}{2\sqrt{\mu(t+2)\frac{L_1^2}{\mu(t+1)}}} \leq 1\ .\nonumber
\end{equation}
Then the previous inequality becomes
\begin{equation}
    F(z_t)-F(z^*) \leq \frac{\mu t+\frac{L_1^2}{\mu(t+1)}}{2}\|z_t-z^*\|^2-\frac{\mu(t+2)+\frac{L_1^2}{\mu(t+1)}}{2}\|z_{t+1}-z^*\|^2+\frac{2L_0^2}{\mu(t+2)}\ .\nonumber
\end{equation}
Multiplying through by $(t+1)$ ensures $(t+1)(F(z_t)-F(z^*))$ is at most
    \begin{equation}
     \frac{\mu t(t+1)+\frac{L_1^2}{\mu}}{2}\|z_t-z^*\|^2-\frac{\mu(t+1)(t+2)+\frac{L_1^2}{\mu}}{2}\|z_{t+1}-z^*\|^2+\frac{2L_0^2}{\mu}\ .\nonumber
\end{equation}
Similarly, from the $\mu$-strongly convex constraint $G(z)$, when $t \in J$, $(t+1)(G(z_t)-G(z^*))$ is at most
\begin{equation}
     \frac{\mu t(t+1)+\frac{L_1^2}{\mu}}{2}\|z_t-z^*\|^2-\frac{\mu(t+1)(t+2)+\frac{L_1^2}{\mu}}{2}\|z_{t+1}-z^*\|^2+\frac{2L_0^2}{\mu}\ .\nonumber
\end{equation}
Summing the two inequalities above up for $t=0,1,2,...,T-1$ yields
\begin{equation}
    \sum_{t \in I}(t+1)(F(z_t)-F(z^*))+\sum_{t \in J}(t+1)(G(z_t)-G(z^*)) \leq \frac{2L_0^2T}{\mu}+\frac{L_1^2\|z_0-z^*\|^2}{2\mu}\ .\nonumber
\end{equation}
For $t \in J$, by definition, we have $G(z_t)>\tau$. Since $G(z^*) \leq 0$, the gap $G(z_t)-G(z^*)>\tau$ is bounded. Then the above inequality becomes
\begin{equation}
    \sum_{t \in I}(t+1)(F(z_t)-F(z^*))+\sum_{t \in J}(t+1)\tau \leq \frac{2L_0^2T}{\mu}+\frac{L_1^2\|z_0-z^*\|^2}{2\mu}.\nonumber
\end{equation}
Therefore, with $T \geq \max\left\{\frac{8L_0^2}{\mu\tau},\sqrt{\frac{2L_1^2\|z_0-z^*\|^2}{\mu\tau}}\right\}$, we have
\begin{align*}
    \sum_{t \in I}&(t+1)(F(z_t)-F(z^*))\\
    &\leq \sum_{t \in I}(t+1)\tau-\sum_{t=0}^{T-1}(t+1)\tau+\frac{2L_0^2T}{\mu}+\frac{L_1^2\|z_0-z^*\|^2}{2\mu} \nonumber\\
    &= \sum_{t \in I}(t+1)\tau-\frac{T(T+1)}{2}\tau+\frac{2L_0^2T}{\mu}+\frac{L_1^2\|z_0-z^*\|^2}{2\mu} \nonumber\\
    &= \sum_{t \in I}(t+1)\tau-\frac{T\tau}{4}\left(T-\frac{8L_0^2}{\mu\tau}\right)-\frac{\tau}{4}\left(T^2-\frac{2L_1^2\|z_0-z^*\|^2}{\mu\tau}\right) \nonumber\\
    &< \sum_{t \in I}(t+1)\tau\ .
\end{align*}
The convexity of $F(z)$ gives us the claimed objective gap bound
\begin{equation}
    F(\Bar{z}_T)-F(z^*)=F\left(\frac{\sum_{t \in I}(t+1)z_t}{\sum_{t \in I}(t+1)}\right)-F(z^*) \leq \frac{\sum_{t \in I}(t+1)F(z_t)}{\sum_{t \in I}(t+1)}-F(z^*)<\tau\ .\nonumber
\end{equation}
The convexity of $G(z)$ gives us the claimed infeasibility bound
\begin{equation}
    G(\Bar{z}_T)=G\left(\frac{\sum_{t \in I}(t+1)z_t}{\sum_{t \in I}(t+1)}\right) \leq \frac{\sum_{t \in I}(t+1)G(z_t)}{\sum_{t \in I}(t+1)}<\tau\ .\nonumber
\end{equation}

\subsection{Proof of Theorem \ref{theorem_mainFJ}}

According to Lemma \ref{lemmafea}, our iterates $x_k$ are always feasible, that is $g(x_k) \leq 0$, for the main problem \eqref{mainproblem} provided $\hat x_{k+1}$ is not an $\epsilon$-FJ point. Note that if $\hat x_{k+1}$ is an $\epsilon$-FJ point, $x_{k}$ must be an $(\epsilon,\epsilon)$-FJ point. For each $x_k$, let $\gamma_{k0}$ and $\gamma_k$ be the necessary multipliers~\eqref{FJoriginal} certifying the optimality of $\hat x_{k+1}$ for the proximal subproblem \eqref{subproblem}. Denote the weighted average of objective and constraint functions for each subproblem as
\begin{equation}
    \mathcal{L}_k(x)=\gamma_{k0}F_k(x)+\gamma_kG_k(x)=\gamma_{k0}\left(f(x)+\frac{\hat{\rho}}{2}\|x-x_k\|^2\right)+\gamma_k\left(g(x)+\frac{\hat{\rho}}{2}\|x-x_k\|^2\right)\ .
    \label{lagrange1}
\end{equation}
Without loss of generality, suppose $\gamma_{k0} \geq 0$, $\gamma_k \geq 0$, and $\gamma_{k0}+\gamma_k=1$. According to FJ conditions~\eqref{FJoriginal}, there exists $\hat{\zeta}_{Fk} \in \partial F_k(\hat{x}_{k+1})$ and $\hat{\zeta}_{Gk} \in \partial G_k(\hat{x}_{k+1})$ which satisfies
\begin{equation}
    \gamma_{k0}\hat{\zeta}_{Fk}+\gamma_k \hat{\zeta}_{Gk} \in -N_X(\hat{x}_{k+1})\ .
    \label{FJcondition1}
\end{equation}
Since $\mathcal{L}_k(x)$ is $(\hat{\rho}-\rho)$-strongly convex, we have
\begin{align*}
    \gamma_{k0}F_k(x_k)+\gamma_kG_k(x_k) \geq &\gamma_{k0}F_k(\hat{x}_{k+1})+\gamma_kG_k(\hat{x}_{k+1}) \\
    &+(\gamma_{k0}\hat{\zeta}_{Fk}+\gamma_k \hat{\zeta}_{Gk})^T(x_k-\hat{x}_{k+1})+\frac{\hat{\rho}-\rho}{2}\|x_k-\hat{x}_{k+1}\|^2\ .
\end{align*}
According to FJ conditions, we also have $\gamma_kG_k(\hat{x}_{k+1})=0$. By \eqref{FJcondition1} and since $x_k \in X$, we know $(\gamma_{k0}\hat{\zeta}_{Fk}+\gamma_k \hat{\zeta}_{Gk})^T(x_k-\hat{x}_{k+1}) \geq 0$. Since $g(x_k) \leq 0$ from Lemma \ref{lemmafea}, the previous inequality becomes
\begin{equation}
    \gamma_{k0}f(x_k) \geq \gamma_{k0}F_k(\hat{x}_{k+1})+\frac{\hat{\rho}-\rho}{2}\|\hat{x}_{k+1}-x_k\|^2\ .\nonumber
\end{equation}
Since $x_{k+1}$ is a $(\tau, \tau)$-solution for the subproblem~\eqref{subproblem}, $F_k(x_{k+1})-F_k(\hat{x}_{k+1}) \leq \tau$. Then the previous inequality becomes
\begin{align}
    \gamma_{k0}f(x_k) &\geq \gamma_{k0}\left(f(x_{k+1})+\frac{\hat{\rho}}{2}\|x_{k+1}-x_k\|^2-\tau\right)+\frac{\hat{\rho}-\rho}{2}\|\hat{x}_{k+1}-x_k\|^2 \nonumber\\
    &\geq \gamma_{k0}(f(x_{k+1})-\tau)+\frac{\hat{\rho}-\rho}{2}\|\hat{x}_{k+1}-x_k\|^2\ .\nonumber
\end{align}
Thus we attain a lower bound for the descent of each step as
\begin{equation}
    \gamma_{k0}(f(x_k)-f(x_{k+1})) \geq \frac{\hat{\rho}-\rho}{2}\|\hat{x}_{k+1}-x_k\|^2-\gamma_{k0}\tau\ . \nonumber
    \label{flborigin}
\end{equation}
When $\gamma_{k0}=0$, then $\|\hat{x}_{k+1}-x_k\|=0$ and so $x_k$ is an exact stationary point of~\eqref{mainproblem}. Now we consider the case that $\gamma_{k0}>0$ here. According to Lemma \ref{lemmastationarity}, before $\hat{x}_{k+1}$ is an $\epsilon$-FJ point, $\|\hat{x}_{k+1}-x_k\|>\frac{\epsilon}{\hat{\rho}}$, then our choice of $\tau=\frac{(\hat{\rho}-\rho)\epsilon^2}{4\hat{\rho}(2\hat{\rho}-\rho)}$ as in \eqref{paraFJ} ensures that
\begin{align}
    f(x_k)-f(x_{k+1}) &\geq \frac{\hat{\rho}-\rho}{2\gamma_{k0}}\|\hat{x}_{k+1}-x_k\|^2-\tau \nonumber\\
    &\geq \frac{\hat{\rho}-\rho}{2}\|\hat{x}_{k+1}-x_k\|^2-\frac{(\hat{\rho}-\rho)\epsilon^2}{4\hat{\rho}(2\hat{\rho}-\rho)} \nonumber\\
    &>\frac{\hat{\rho}-\rho}{2} \frac{\epsilon^2}{\hat{\rho}^2}-\frac{(\hat{\rho}-\rho)\epsilon^2}{4\hat{\rho}(2\hat{\rho}-\rho)} \nonumber\\
    &>\frac{(\hat{\rho}-\rho)\epsilon^2}{2\hat{\rho}^2}-\frac{(\hat{\rho}-\rho)\epsilon^2}{4\hat{\rho}^2} \nonumber\\
    &=\frac{(\hat{\rho}-\rho)\epsilon^2}{4\hat{\rho}^2}\ .\nonumber
\end{align}
Hence by Assumption \ref{assumption2}, the number of total iterations $K$ of Algorithm \ref{algorithm_main} before an $(\epsilon,\epsilon)$-FJ point is found is upper bounded by
\begin{equation}
    K<\frac{4\hat{\rho}^2(f(x_0)-f_{lb})}{(\hat{\rho}-\rho)\epsilon^2}\ .\nonumber
\end{equation}
Consequently, Algorithm~\ref{algorithm_main} (which uses Algorithm~\ref{algorithm_sub} for $T$ steps as an oracle each iteration) will identify an $(\epsilon,\epsilon)$-FJ point using at most the following total number of subgradient evaluations of either the objective or constraints
\begin{equation}
    KT<\frac{4\hat{\rho}^2(f(x_0)-f_{lb})}{(\hat{\rho}-\rho)\epsilon^2}\max\left\{\frac{96\hat{\rho}(2\hat{\rho}-\rho)(3M^2-2\hat{\rho}g_{lb})}{(\hat{\rho}-\rho)^2\epsilon^2},\sqrt{\frac{288\hat{\rho}^3(2\hat{\rho}-\rho)D^2}{(\hat{\rho}-\rho)^2\epsilon^2}}\right\} \ .\nonumber
\end{equation}

\subsection{Proof of Theorem \ref{theorem_mainKKT}}
Nearly identical reasoning to that of Theorem \ref{theorem_mainFJ}'s proof under the constraint qualification Assumption~\ref{assumption5} yields our claimed result of approximate KKT stationarity convergence rate. The exact details for this symmetric case are provided in the appendix for completeness.

\subsection{Proof of Lemmas}

\subsubsection{Proof of Lemma \ref{lemmastationarity}} 
First, we consider the claimed result of approximate Fritz-John stationarity on the original nonsmooth nonconvex problem~\eqref{mainproblem}. Necessarily the FJ conditions~\eqref{FJoriginal} are satisfied for proximally subproblem~\eqref{subproblem} at $\hat{x}_{k+1}$ for some $\gamma_{k0},\gamma_k \geq 0$, $\gamma_{k0}+\gamma_k=1$, $\hat{\zeta}_{Fk} \in \partial F_k(\hat{x}_{k+1})$ and $\hat{\zeta}_{Gk} \in \partial G_k(\hat{x}_{k+1})$. By the sum rule of subgradient calculus, let $\hat{\zeta}_{fk}=\hat{\zeta}_{Fk}-\hat{\rho}(\hat{x}_{k+1}-x_k) \in \partial f(\hat{x}_{k+1})$ and $\hat{\zeta}_{gk}=\hat{\zeta}_{Gk}-\hat{\rho}(\hat{x}_{k+1}-x_k) \in \partial g(\hat{x}_{k+1})$. The FJ conditions for the proximal subproblem guarantee there exists some $\nu \in N_X(\hat{x}_{k+1})$ such that
\begin{align*}
    \gamma_{k0}\left[\hat{\zeta}_{fk}+\hat{\rho}(\hat{x}_{k+1}-x_k)\right]+\gamma_k\left[\hat{\zeta}_{gk}+\hat{\rho}(\hat{x}_{k+1}-x_k)\right]=-\nu \ .
\end{align*}
Hence when $\|\hat{x}_{k+1}-x_k\| \leq \frac{\epsilon}{\hat{\rho}}$, the first approximate FJ condition~\eqref{FJ1_1} holds at $\hat x_{k+1}$ for the original nonsmooth nonconvex problem as $\|\gamma_{k0}\hat{\zeta}_{fk}+\gamma_k\hat{\zeta}_{gk}+\nu\|=\hat{\rho}\|\hat{x}_{k+1}-x_k\| \leq \epsilon.$

Moreover, we can verify the second approximate FJ condition~\eqref{FJ1_2} at $\hat x_{k+1}$ in the following two cases: When $\gamma_k=0$, this trivially holds as $|\gamma_k g(\hat{x}_{k+1})|=0$. When $\gamma_k>0$, we have $G_k(\hat{x}_{k+1})=0$ according to FJ conditions. Hence $0 \geq g(\hat{x}_{k+1})=-\frac{\hat{\rho}}{2}\|\hat{x}_{k+1}-x_k\|^2 \geq -\frac{\epsilon^2}{2\hat{\rho}}.$
As a result,
\begin{equation}
    |\gamma_kg(\hat{x}_{k+1})| \leq |g(\hat{x}_{k+1})| \leq \frac{\epsilon^2}{2\hat{\rho}}<\epsilon^2 \ .\nonumber
\end{equation}

Nearly identical reasoning under the constraint qualification Assumption~\ref{assumption5} yields our claimed result of approximate KKT stationarity on the original nonsmooth nonconvex problem~\eqref{mainproblem}. The exact details for this symmetric case are provided in the appendix for completeness.

\subsubsection{Proof of Lemma \ref{lemma_quadraticbound}} 
Let $z_0=x_k$ and $z^*=\hat{x}_{k+1}$ be the optimal solution for problem \eqref{subproblem}, which is $\mu=\hat{\rho}-\rho$-strongly convex. Consider any $\zeta_{Fk} \in \partial F_k(z)$, $\zeta_{Gk} \in \partial G_k(z)$, which the sum rule ensures have $\zeta_{f}=\zeta_{Fk}-\hat{\rho}(z-z_0) \in \partial f(z)$, and $\zeta_{g}=\zeta_{Gk}-\hat{\rho}(z-z_0) \in \partial g(z)$. 

First, we verify the non-Lipschitz subgradient bound for $F_k$ with $L^2_0\geq 9M^2$ and $L_1 \geq 6\hat{\rho}$. Namely, consider any $z$ with $G_k(z)\leq \tau$. Then
\begin{align*}
    & L_0^2+L_1(F_k(z)-F_k(z^*)) \\ &=L_0^2+L_1(F_k(z)-F_k(z^*))-\|\zeta_{Fk}\|^2+\|\zeta_{Fk}\|^2 \\
    &\geq L_0^2+L_1\left(f(z)+\frac{\hat{\rho}}{2}\|z-z_0\|^2-F_k(z_0)\right)-\|\zeta_{f}+\hat{\rho}(z-z_0)\|^2+\|\zeta_{Fk}\|^2 \\
    &= L_0^2+L_1(f(z)-f(z_0))+\frac{L_1\hat{\rho}}{2}\|z-z_0\|^2-\|\zeta_{f}\|^2-2\hat{\rho}\zeta_{f}^T(z-z_0)-\hat{\rho}^2\|z-z_0\|^2+\|\zeta_{Fk}\|^2 \\
    &\geq L_0^2-L_1M\|z-z_0\|+\frac{L_1\hat{\rho}}{2}\|z-z_0\|^2-M^2-2\hat{\rho}M\|z-z_0\|-\hat{\rho}^2\|z-z_0\|^2+\|\zeta_{Fk}\|^2 \\
    &=(L_0^2-M^2)-(L_1+2\hat{\rho})M\|z-z_0\|+\left(\frac{L_1}{2}-\hat{\rho}\right)\hat{\rho}\|z-z_0\|^2+\|\zeta_{Fk}\|^2 \\
    &\geq \left(L_0^2-M^2-\frac{(L_1+2\hat{\rho})^2M^2}{2(L_1-2\hat{\rho})\hat{\rho}}\right)+\|\zeta_{Fk}\|^2\\
    &\geq \|\zeta_{Fk}\|^2
\end{align*}
where the first inequality uses that $F_k(z_0) \geq F_k(z^*)$, the second uses the $M$-Lipschitz continuity of $f$, the third minimizes over all $\|z-z_0\|$ (noting that $L_1 >2\hat\rho$), and the last inequality uses the assumed bounds on $L_0^2$ and $L_1$.

Similarly, we verify the non-Lipschitz subgradient bound for the proximally penalized constraints $G_k$ with $L_0^2 \geq 9M^2-6\hat{\rho}g(z_0), L_1 \geq 6\hat{\rho}$. Namely, for any $z$ with $G_k(z)>\tau$, we have
\begin{align*}
    &L_0^2+L_1(G_k(z)-G_k(z^*)) \\
    &=L_0^2+L_1(G_k(z)-G_k(z^*))-\|\zeta_{Gk}\|^2+\|\zeta_{Gk}\|^2 \\
    &\geq L_0^2+L_1\left(g(z)+\frac{\hat{\rho}}{2}\|z-z_0\|^2\right)-\|\zeta_{g}+\hat{\rho}(z-z_0)\|^2+\|\zeta_{Gk}\|^2 \\
    &=L_0^2+L_1g(z_0)+L_1(g(z)-g(z_0))+\frac{L_1\hat{\rho}}{2}\|z-z_0\|^2-\|\zeta_{g}\|^2-2\hat{\rho}\zeta_{g}^T(z-z_0)-\hat{\rho}^2\|z-z_0\|^2+\|\zeta_{Gk}\|^2 \\
    &\geq L_0^2+L_1g(z_0)-L_1M\|z-z_0\|+\frac{L_1\hat{\rho}}{2}\|z-z_0\|^2-M^2-2\hat{\rho}M\|z-z_0\|-\hat{\rho}^2\|z-z_0\|^2+\|\zeta_{Gk}\|^2 \\
    &=(L_0^2-M^2+L_1g(z_0))-(L_1+2\hat{\rho})M\|z-z_0\|+\left(\frac{L_1}{2}-\hat{\rho}\right)\hat{\rho}\|z-z_0\|^2+\|\zeta_{Gk}\|^2 \\
    &\geq \left(L_0^2-M^2+L_1g(z_0)-\frac{(L_1+2\hat{\rho})^2M^2}{2(L_1-2\hat{\rho})\hat{\rho}}\right) +\|\zeta_{Gk}\|^2\\
    & \geq \|\zeta_{Gk}\|^2.
\end{align*}
Since $g(z_0) \geq g_{lb}$, setting $L_0^2=9M^2-6\hat{\rho}g_{lb}$ and $L_1=6\hat{\rho}$ satisfies both cases above.

\subsubsection{Proof of Lemma \ref{lemmadualvar}} 
Let $\hat{x}_{k+1}$ be the exact solution for problem \eqref{subproblem} with optimal dual variable $\lambda_k$. The $(\hat{\rho}-\rho)$-strong convexity of $G_k(x)$ implies that the set $\{x|G_k(x) \leq 0\}$ has diameter $D=\sqrt{(-8g_{lb})/(\hat{\rho}-\rho)}$. Since $x_k$ and $\hat{x}_{k+1}$ both lying in this set, $\|\hat{x}_{k+1}-x_k\| \leq D$.
According to KKT conditions \eqref{KKToriginal}, there exists $\hat{\zeta}_{Fk} \in \partial F_k(\hat{x}_{k+1})$ and $\hat{\zeta}_{Gk} \in \partial G_k(\hat{x}_{k+1})$ which satisfies $\hat{\zeta}_{Fk}+\lambda_k \hat{\zeta}_{Gk} \in -N_X(\hat{x}_{k+1})$. Trivially if $\lambda_k$ is zero, it is bounded, so we consider $\lambda_k$ is positive (and so $G_k(\hat{x}_{k+1})=0$). Then there exists $\nu \in N_X(\hat{x}_{k+1})$ such that $\hat{\zeta}_{Fk}+\lambda_k\hat{\zeta}_{Gk} =-\nu$. Hence
\begin{align}
    \lambda_k&=\frac{\|\hat{\zeta}_{Fk}\|}{\|\hat{\zeta}_{Gk}+\frac{\nu}{\lambda_k}\|}.
    \label{lambda1}
\end{align}
We can directly upper bound the numerator above as $\hat{\zeta}_{fk}=\hat{\zeta}_{Fk}-\hat{\rho}(\hat{x}_{k+1}-x_k) \in \partial f_k(\hat{x}_{k+1})$. So Assumption \ref{assumption3} and the bound $\|\hat{x}_{k+1}-x_k\| \leq D$, ensure $\|\hat{\zeta}_{Fk}\| \leq M+\hat\rho D$.
Assumption \ref{assumption5} facilitates lower bounding the denominator above. Namely, there must exist $v \in -N_X^*(\hat{x}_{k+1})$ with $\|v\|=1$, such that $\hat{\zeta}_{Gk}^Tv \leq -\sigma$. Since $\nu \in N_X(\hat{x}_{k+1})$ and $v \in -N_X^*(\hat{x}_{k+1})$, we know $\nu^Tv \leq 0$. Then
$\|\hat{\zeta}_{Gk}+\frac{\nu}{\lambda_k}\|=\|\hat{\zeta}_{Gk}+\frac{\nu}{\lambda_k}\| \cdot \|v\| \geq -(\hat{\zeta}_{Gk}+\frac{\nu}{\lambda_k})^Tv \geq \sigma.$ Combining these upper and lower bounds gives the claimed uniform Lagrange multiplier bound.

\subsubsection{Proof of Lemma \ref{lemmafea}} 
    
First, we inductively show the feasibility of the iterates $x_k$ before $\hat x_{k+1}$ is an $\epsilon$-FJ point. Assume $G_k(x_k)=g(x_k) \leq 0$. Necessarily the FJ conditions~\eqref{FJoriginal} are satisfied for proximally subproblem~\eqref{subproblem} at $\hat{x}_{k+1}$ for some $\gamma_{k0},\gamma_k \geq 0$, $\gamma_{k0}+\gamma_k=1$, $\hat{\zeta}_{Fk} \in \partial F_k(\hat{x}_{k+1})$ and $\hat{\zeta}_{Gk} \in \partial G_k(\hat{x}_{k+1})$. Consider the function $\mathcal{L}_k(x) = \gamma_{k0}F_k(x)+\gamma_kG_k(x)$, which is minimized over $X$ at $\hat x_{k+1}$. The $(\hat{\rho}-\rho)$-strong convexity of $F_k$ and $G_k$ ensures
\begin{align*}
    \gamma_{k0}F_k(x_{k+1})+\gamma_kG_k(x_{k+1})
    \geq &\gamma_{k0}F_k(\hat{x}_{k+1})+\gamma_kG_k(\hat{x}_{k+1})+(\gamma_{k0}\hat{\zeta}_{Fk} \\
    &+\gamma_k \hat{\zeta}_{Gk})^T(x_{k+1}-\hat{x}_{k+1})+\frac{\hat{\rho}-\rho}{2}\|x_{k+1}-\hat{x}_{k+1}\|^2.
\end{align*}
The Fritz-John conditions ensure that $\gamma_k G_k(\hat{x}_{k+1})=0$ and that $\lambda_{k0}\hat{\zeta}_{Fk}+\lambda_k \hat{\zeta}_{Gk} \in -N_X(\hat{x}_{k+1})$, which guarantees $x_{k+1} \in X$ has $(\gamma_{k0}\zeta_{Fk}+\gamma_k \zeta_{Gk})^T(x_{k+1}-\hat{x}_{k+1}) \geq 0$. These two observations simplify the above inequality to
\begin{align*}
    \gamma_{k0}(F_k(x_{k+1}) - F_k(\hat{x}_{k+1}))+\gamma_kG_k(x_{k+1})
    \geq \frac{\hat{\rho}-\rho}{2}\|x_{k+1}-\hat{x}_{k+1}\|^2 \ .
\end{align*}
By Corollary~\ref{corollary_sub}, the proposed selection of $T_{inner}$ ensures $x_{k+1}$ is a $(\tau,\tau)$-optimal solution for the subproblem \eqref{subproblem}. Hence $F_k(x_{k+1})-F_k(\hat{x}_{k+1}) \leq \tau$ and $G_k(x_{k+1}) \leq \tau$. Noting $\gamma_{k0}+\gamma_k=1$, the above inequality further simplifies to
\begin{align}
    \|\hat{x}_{k+1}-x_{k+1}\| &\leq \sqrt{\frac{2\tau}{\hat{\rho}-\rho}} \ .\nonumber
\end{align}
Assume $\hat x_{k+1}$ is not an $\epsilon$-FJ point, according to Lemma~\ref{lemmastationarity}, $\|\hat{x}_{k+1}-x_k\|>\frac{\epsilon}{\hat{\rho}}$. Thus
\begin{equation}
    \|x_{k+1}-x_k\|^2 \geq \frac{1}{2}\|\hat{x}_{k+1}-x_k\|^2-\|\hat{x}_{k+1}-x_{k+1}\|^2>\frac{\epsilon^2}{2\hat{\rho}^2}-\frac{2\tau}{\hat{\rho}-\rho} \ .\nonumber
\end{equation}
By our selection of $\tau=\frac{(\hat{\rho}-\rho)\epsilon^2}{4\hat{\rho}(2\hat{\rho}-\rho)}$ as in \eqref{paraFJ}, every iteration prior to finding an  $\epsilon$-FJ point must have
\begin{equation}
    \|x_{k+1}-x_k\|^2>\frac{(\hat{\rho}-\rho)\epsilon^2}{2\hat{\rho}^2(2\hat{\rho}-\rho)} \ .
    \label{stoppingcriteriaFJ2}
\end{equation}
Therefore $g(x_{k+1}) \leq 0$ is inductively ensured if $g(x_k) \leq 0$ and $\hat x_{k+1}$ is not an $\epsilon$-FJ point as
\begin{align*}
    g(x_{k+1})&=G(x_{k+1})-\frac{\hat{\rho}}{2}\|x_{k+1}-x_k\|^2 
    \leq \tau-\frac{\hat{\rho}}{2} \frac{(\hat{\rho}-\rho)\epsilon^2}{2\hat{\rho}^2(2\hat{\rho}-\rho)} 
    =0 \ .
\end{align*}

By nearly identical reasoning, we find that under the KKT parameter selections~\eqref{paraKKT}, the feasibility of $x_k$ ensures $x_{k+1}$ is feasible so long as $\hat x_{k+1}$ is not an $\epsilon$-KKT point. The details of this are deferred to the appendix for completeness.

    \section{Numerics with Sparsity Inducing SCAD Constraints} \label{numerical}
Lastly, we illustrate the diversity of approximate stationary points actually reached by the inexact proximal point method. {\color{\blk} Note that our method is not designed to be state-of-the-art in terms of actual computational performance; rather, here, we aim to showcase its diversity of possible limiting behaviors, all of which our theory captures.} The frequent occurrences of FJ points (numerically failing to have MFCQ) seen here support our work and motivate future works developing methods capable of handling such limit points {\color{black} in nonsmooth optimization settings}.

We consider the sparse phase retrieval (SPR) problem previously described in Section \ref{subsec:vignette}. Phase retrieval is a common problem in various applications, such as imaging, X-ray crystallography, and transmission electron microscopy. The phase is recovered by solving linear equations $Ax=b$ up to sign changes, $(Ax)^2=b^2$. We construct our sparse phase retrieval problem as
\begin{equation}
\begin{cases}
    \min_{x \in X} & f(x)=\frac{1}{m}\sum_{i=1}^{m} |(a_i^Tx)^2-b_i^2| \\
    \mathrm{s.t.} & g(x)=\sum_{i=1}^{n} s(x_i)-p \leq 0\ .
    \label{SPR}
\end{cases}
\end{equation}
Here $s: \mathbb{R} \to \mathbb{R}$ is the Smoothly Clipped Absolute Deviation (SCAD) function~\eqref{scad}.
{\color{black} Although the objective and the constraints are nonsmooth, this only occurs on a measure zero set. By starting our experiments from a random initialization, we ensure such points are never reached. Consequently, the subgradient oracle we use is exactly the gradient.} Despite the simple piecewise quadratic definition of SCAD constraints, whenever $p$ is a multiple of three, proximal subproblems exist where MFCQ fails (that is, no Slater points exist). Consider $x_k=[5,5,\cdots,5,0,0,\cdots,0]$ which consists of $p/3$ fives and $(n-p/3)$ zeroes. Then the MFCQ condition fails as $G_k(x_k)=0$ as $0 \in \partial G_k(x_k)$.

In Section~\ref{subsec:stopping}, we first discuss our synthetic SPR problem instances and propose a simple stopping criterion, which we find is numerically effective. Then Section~\ref{subsec:experiments} presents numerical results from applying our Proximally Guided Switching Subgradient Method to SPR problems, identifying varied convergence to FJ points, KKT points with active constraints, and KKT points with inactive constraints.

\subsection{SPR Problem Generation, Stopping Criteria and Experimental Evaluation} \label{subsec:stopping}
Note the SPR problem~\eqref{SPR} has $f: \mathbb{R}^{n} \to \mathbb{R}$ and $g: \mathbb{R}^{n} \to \mathbb{R}$ weakly convex and Lipschitz on $X=[-10,10]^n$. For our numerics, we fix $m=120$, $n=120$, and $p \in [0, 3n)$ varies to control the sparsity of our problem. Then $A \in \mathbb{R}^{m \times n}$ has entries sampled independently from a standard Gaussian distribution. To have a sparse optimal solution, sample values for the first $30$ indices of $x^*$ uniformly in $\pm [5,10]$, and set the other $90$ entries as $0$. We then compute $b^2=(Ax^*)^2+\eta$ with noise $\eta$ drawn from a standard Gaussian. Our numerics use a random feasible initialization $x_0$ with entries sampled from $N(0,0.01)$ independently.

According to Lemma B.1 in \cite{davis2019proximally}, $f(x)$ is expected to be $2$-weakly convex. To leave some gap, we set $\rho=3$, and $\hat{\rho}=6$. As other inputs to the Proximally Guided Switching Subgradient Method, we set $\epsilon=0.01$ and run the method for $K=10^3$ outer iterations, each using $T=10^4$ inner steps. Consequently, we use a total of $10^7$ subgradient evaluations.

\subsubsection{Stopping Criteria}
Our convergence theory supporting Theorems~\ref{theorem_mainFJ} and~\ref{theorem_mainKKT} showed that the iterates $x_k$ of the Proximally Guided Switching Subgradient Method stay feasible and are guaranteed to decrease the objective value until an $(\epsilon,\epsilon)$-FJ or KKT point is found. This motivates the following simple stopping criterion: continue taking inexact proximal steps until either
\begin{equation}
   \quad g(x_k)>0 \quad \text{or} \quad f(x_k) \geq f(x_{k-1}) \ . \nonumber
    \label{stoppingcriterianumerics}
\end{equation}
In the following numerics, we denote the first time this condition is reached via a vertical dotted line. We find numerically that this criterion aligns well with when the associated Fritz-John and Lagrange multipliers and the iterate's feasibility level out.

This stopping criterion is never satisfied prior to reaching an approximate stationary point (see Lemma~\ref{lemmafea}), and so Algorithm \ref{algorithm_main} continues. The first time the stopping criterion is met, we must have reached our targeted approximate stationary point and stopped our algorithm. Generally, however, the stopping criterion may fail to be satisfied despite the iterates being an $(\epsilon,\epsilon)$-FJ or KKT point, so these conditions are heuristic in nature.

{\color{black} \subsubsection{Experimental Evaluation}
For evaluation purposes, approximate Lagrange multipliers can be directly extracted from the subproblem solves. Numerically, we found iterates almost never reached the boundary of $X=[-10,10]^n$ and hence we omit projection steps from our discussion for ease. For $t\in \mathbb{N}$, let $z_t$, $\zeta_t$ and $\alpha_t$ denote the iterates, subgradients, and step-sizes in iteration $k$'s call to the switching subgradient method with $\tau=0$. Given $\sum_t \alpha_t$ diverges and $z_t$ converges to the unique minimizer $\hat x_{k+1}$, letting $I_F=\{t \in \mathbb{N} | G(z_t)\leq 0\}$, 
$$ 0 = \frac{x_k - \hat x_{k+1} }{\sum_t \alpha_t}  = \frac{\sum_{t \in I_F}\alpha_t\zeta_{Ft} + \sum_{t \notin I_F}\alpha_t \zeta_{Gt}}{\sum_t \alpha_t} \ . $$

Plugging in that $\zeta_{Ft} = \zeta_{ft} + \hat{\rho}(z_t - x_k)$ and $\zeta_{Gt} = \zeta_{gt} + \hat{\rho}(z_t - x_k)$ and again using $z_t$ converging towards $\hat x_{k+1}$, we conclude
$$ \frac{\sum_{t \in I_F}\alpha_t\zeta_{ft} + \sum_{t \notin I_F}\alpha_t \zeta_{gt}}{\sum_t \alpha_t} = \hat\rho (x_k - \hat x_{k+1}) $$

As $x_k - \hat x_{k+1}$ tends to zero (i.e., as the outer proximal point iteration converges), the lefthand-side converging to zero is approximately the needed condition for FJ stationarity. Consequently, natural choices of Fritz-John and Lagrange multipliers from iteration $k$'s inner loop are then
$$\gamma_{k0} = \frac{\sum_{t \in I_F}\alpha_t}{\sum_t \alpha_t}, \quad \gamma_k = \frac{\sum_{t \notin I_F}\alpha_t}{\sum_t \alpha_t}, \quad \lambda_k = \frac{\sum_{t \notin I_F}\alpha_t}{\sum_{t \in I_F} \alpha_t} \ . $$
As the switching subgradient method's iterates converge, $\zeta_{fk} = \sum_{t \in I_F}\alpha_t\zeta_{ft}\big/\left(\sum_{t \in I_F} \alpha_t\right)$ and $\zeta_{gk} = \sum_{t \not\in I_F}\alpha_t\zeta_{gt}\big/\big(\sum_{t \not\in I_F} \alpha_t\big)$ converge to elements of $\partial f(\hat x_{k+1})$ and $\partial g(\hat x_{k+1})$. Using these multipliers and limiting subgradients, the associated FJ and KKT stationarity measures for $f$ and $g$ are then $\|\gamma_{k0}\zeta_{fk} + \gamma_k\zeta_{gk}\| = \hat{\rho}\|\hat x_{k+1}-x_k\|$ and $\|\zeta_{fk} + \lambda_k\zeta_{gk}\| = (1+\lambda_k)\hat{\rho}\|\hat x_{k+1}-x_k\|$. For our experiments, we approximate the above formula for multipliers by using $T$ iterations of the switching subgradient method with $\tau =O(\epsilon^2)$ (as in~\eqref{paraFJ} and~\eqref{paraKKT}) and approximate the above stationarity measures using $x_{k+1}$ as a proxy for $\hat x_{k+1}$.}

\subsection{Three Distinct Families of Limit Points in Sparse Phase Retrieval} \label{subsec:experiments}
For our randomly generated SPR problem instances, we consider three different selections of $p$, namely $90,91,320$. Although our iterates always converge to an approximate FJ or KKT point, we can see three distinct behaviors under different levels of sparsity controlled by the SCAD constraint. {\color{black} When $p$ is small $p=90$ and $p=91$, we see convergence to a range of FJ and KKT limit points on the boundary of the constraint set. When $p$ is large, we often see convergence to the interior of the feasible region (with the constraint being inactive). Consequently, the Lagrange multiplier tends to zero}, so approximate FJ and KKT stationarities are equivalent. The MATLAB source code implementing these experiments is available at \url{https://github.com/Zhichao-Jia/arXiv_proximal2022}.

\begin{figure}[htbp]
    \centering
    \begin{subfigure}{0.3\textwidth}
    {
        \centering
        \includegraphics[width=\textwidth]{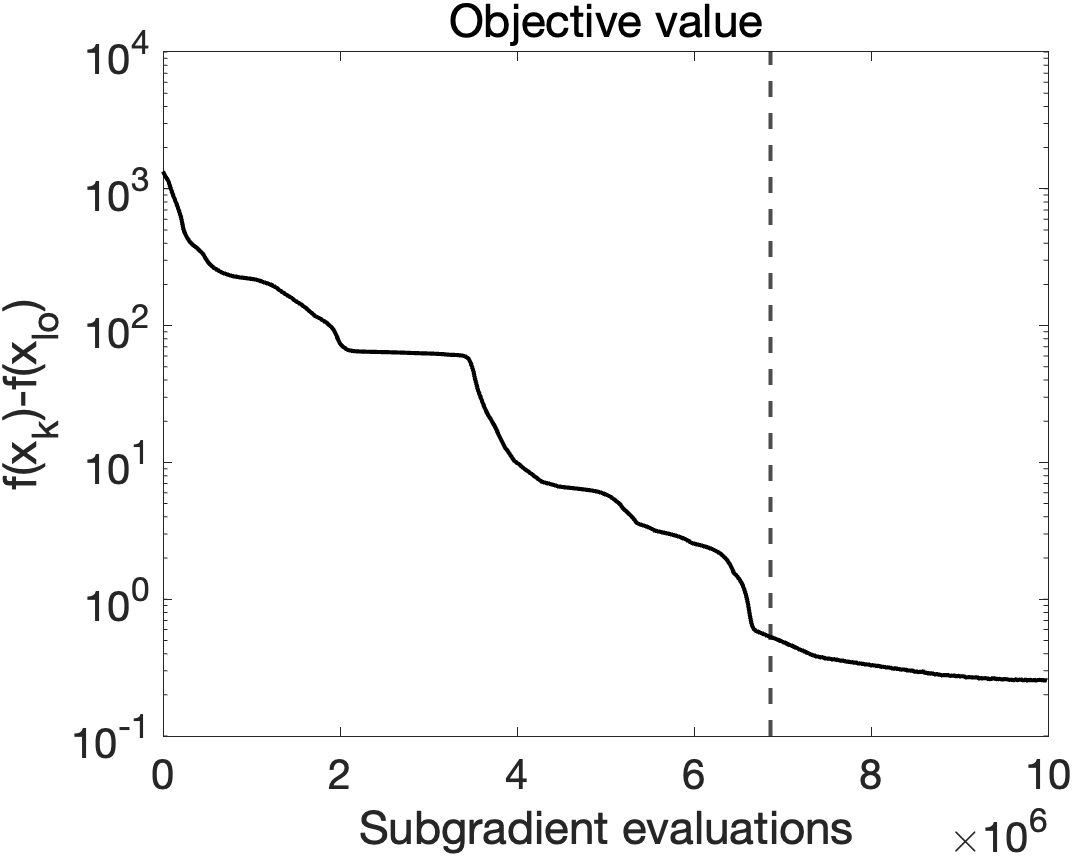}
        \caption{}
        \label{spr1a}
    }
    \end{subfigure}
    \hfill
    \begin{subfigure}{0.3\textwidth}
    {
        \centering
        \includegraphics[width=\textwidth]{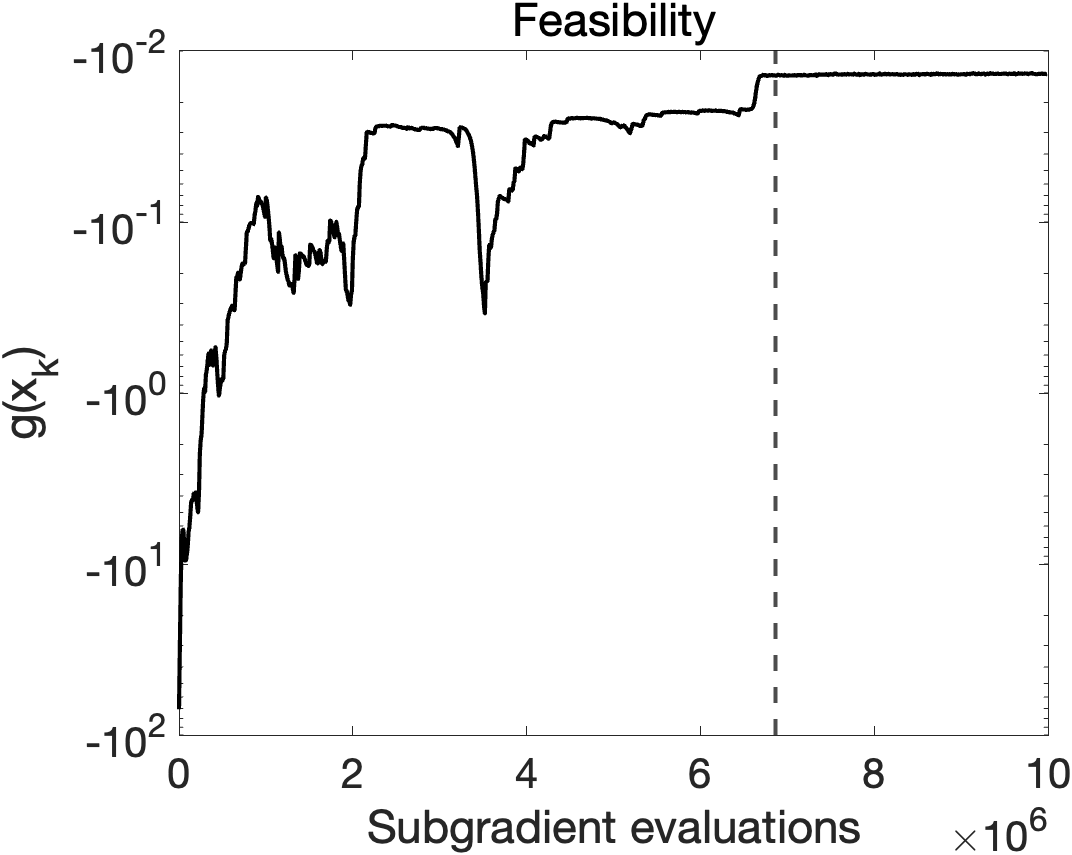}
        \caption{}
        \label{spr1b}
    }
    \end{subfigure}
    \hfill
    \begin{subfigure}{0.3\textwidth}
    {
        \centering
        \includegraphics[width=\textwidth]{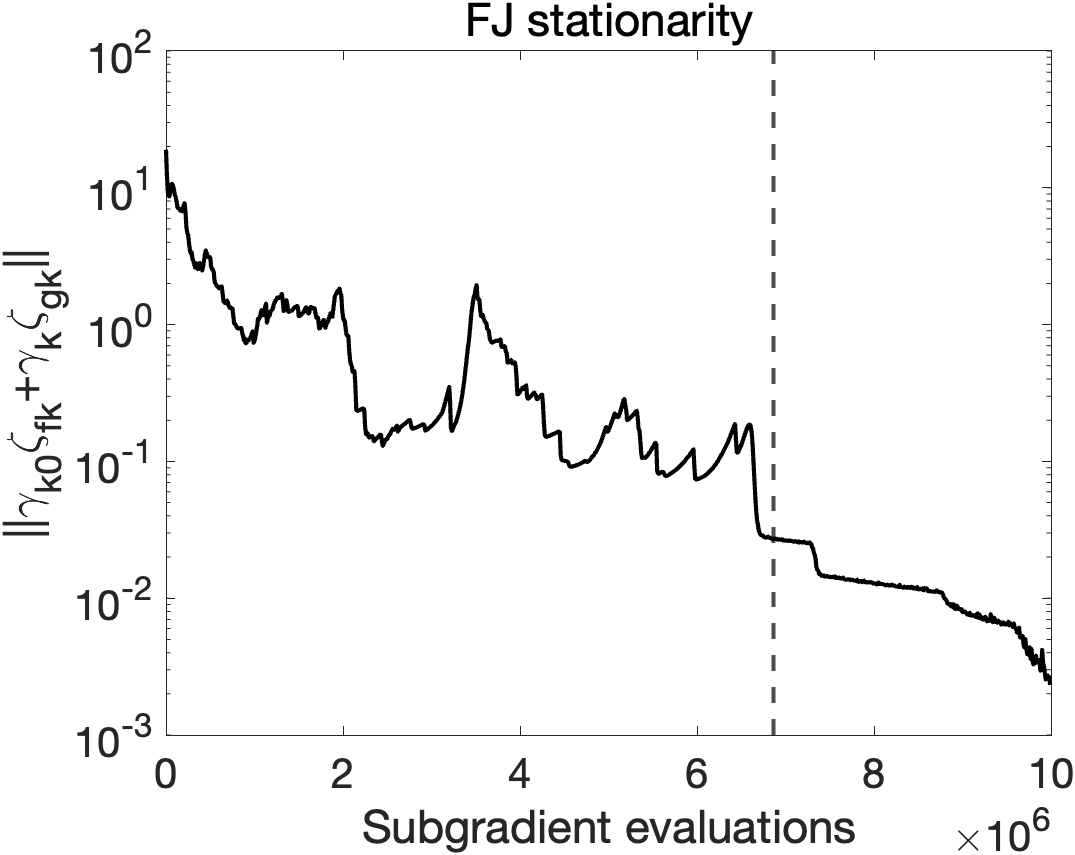}
        \caption{}
        \label{spr1c}
    }
    \end{subfigure}
    \begin{subfigure}{0.3\textwidth}
    {
        \centering
        \includegraphics[width=\textwidth]{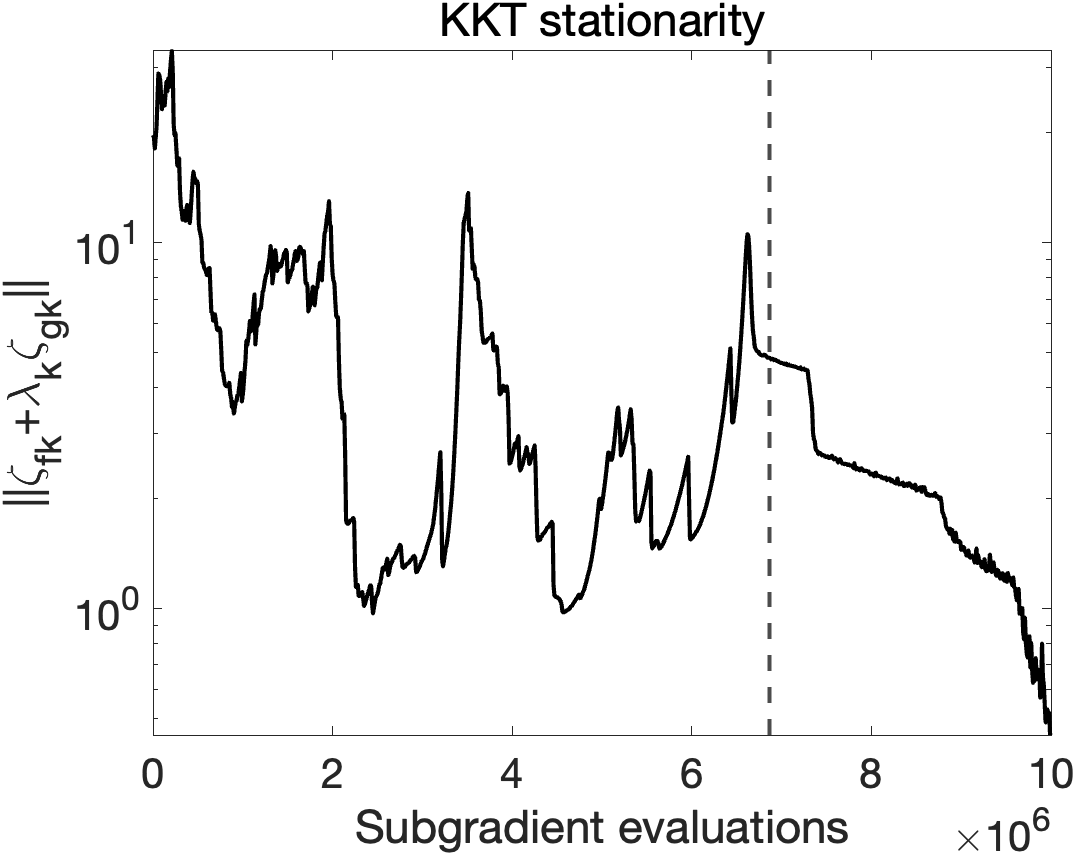}
        \caption{}
        \label{spr1d}
    }
    \end{subfigure}
    \hfill
    \begin{subfigure}{0.3\textwidth}
    {
        \centering
        \includegraphics[width=\textwidth]{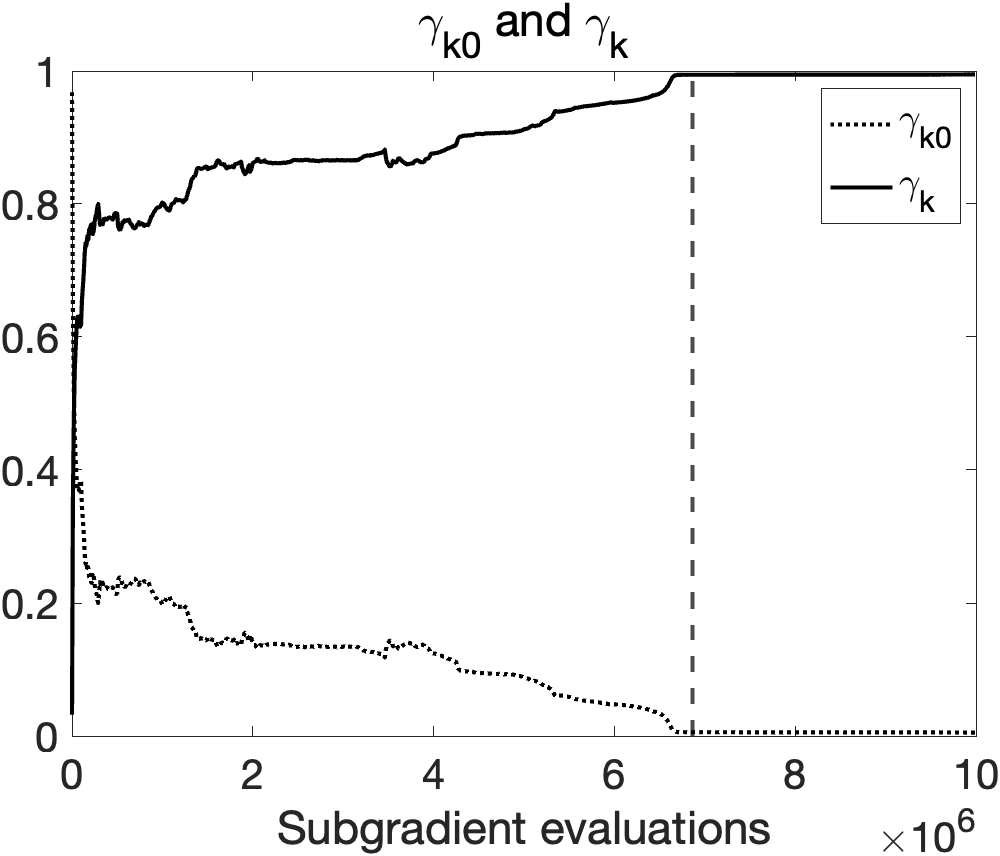}
        \caption{}
        \label{spr1e}
    }
    \end{subfigure}
    \hfill
    \begin{subfigure}{0.3\textwidth}
    {
        \centering
        \includegraphics[width=\textwidth]{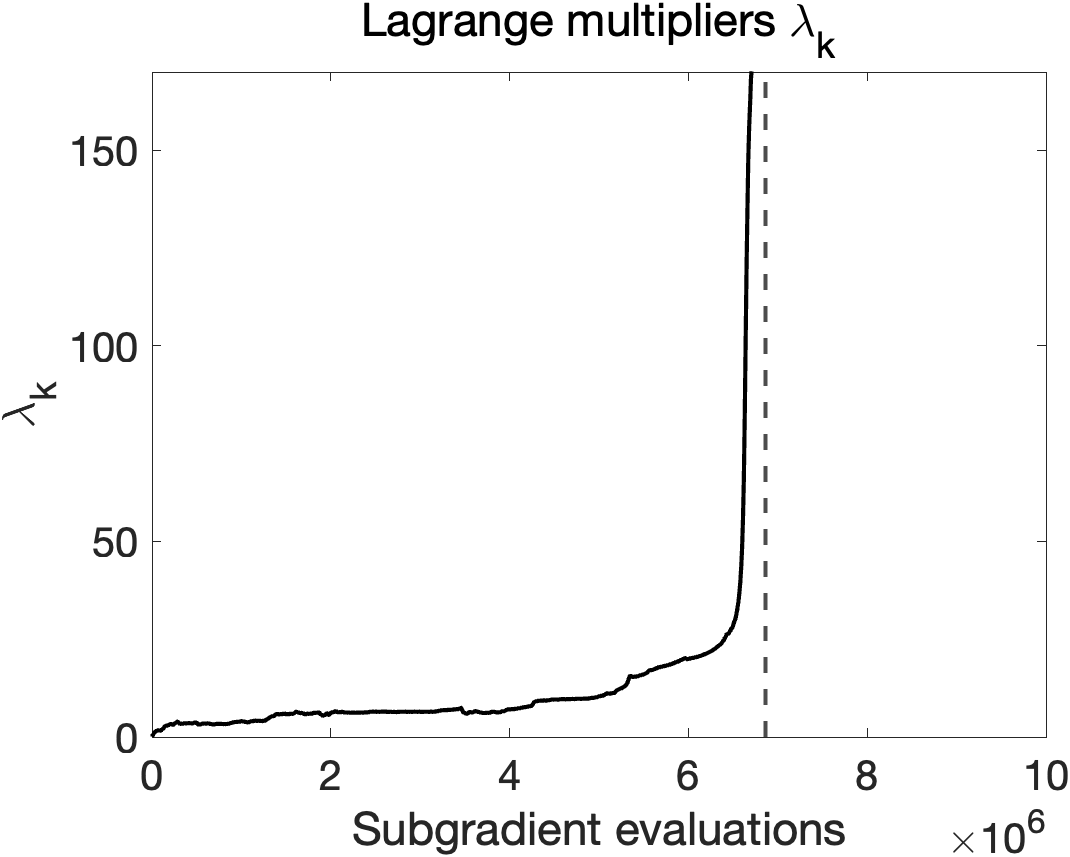}
        \caption{}
        \label{spr1f}
    }
    \end{subfigure}
    \caption{Finding FJ Stationarity: $K=10^3,T=10^4,p=90$. Dotted lines show where the stopping criteria applied. $x_{lo}$ is the stationary point near the final iterate.}
    \label{example1}
\end{figure}

\begin{figure}[htbp]
    \centering
    \begin{subfigure}{0.3\textwidth}
    {
        \centering
        \includegraphics[width=\textwidth]{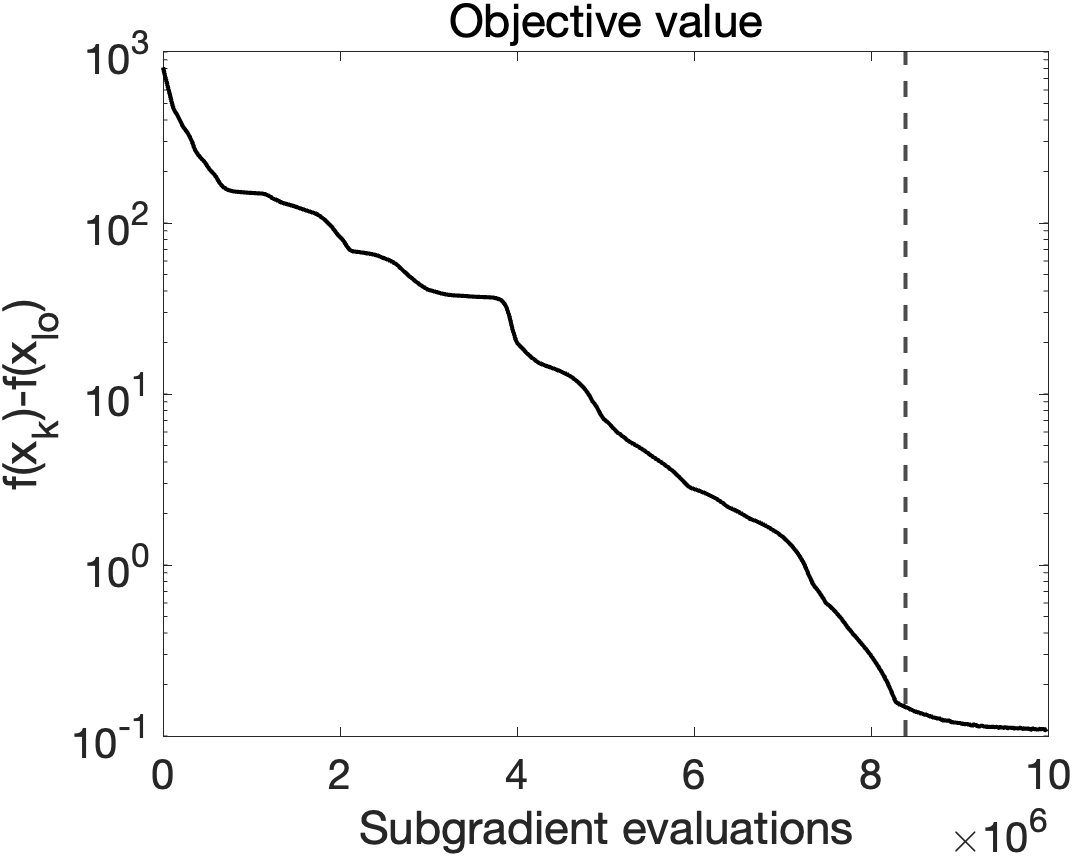}
        \caption{}
        \label{spr2a}
    }
    \end{subfigure}
    \hfill
    \begin{subfigure}{0.3\textwidth}
    {
        \centering
        \includegraphics[width=\textwidth]{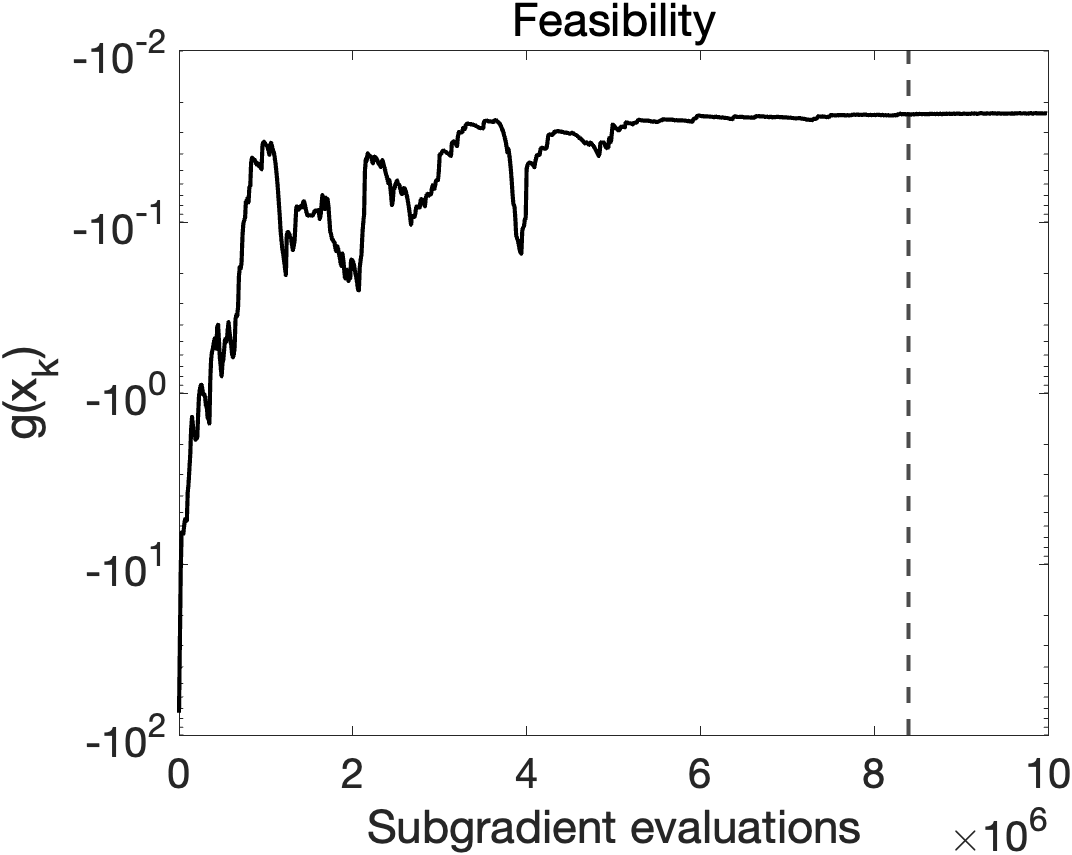}
        \caption{}
        \label{spr2b}
    }
    \end{subfigure}
    \hfill
    \begin{subfigure}{0.3\textwidth}
    {
        \centering
        \includegraphics[width=\textwidth]{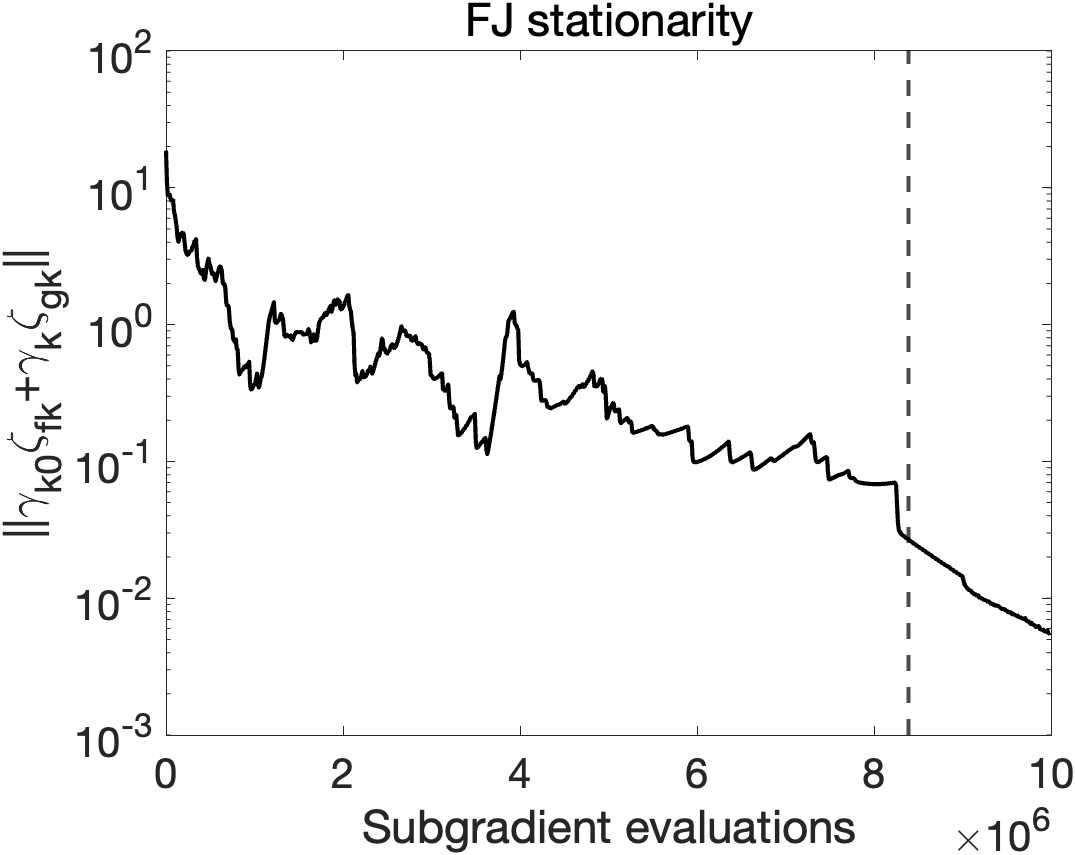}
        \caption{}
        \label{spr2c}
    }
    \end{subfigure}
    \begin{subfigure}{0.3\textwidth}
    {
        \centering
        \includegraphics[width=\textwidth]{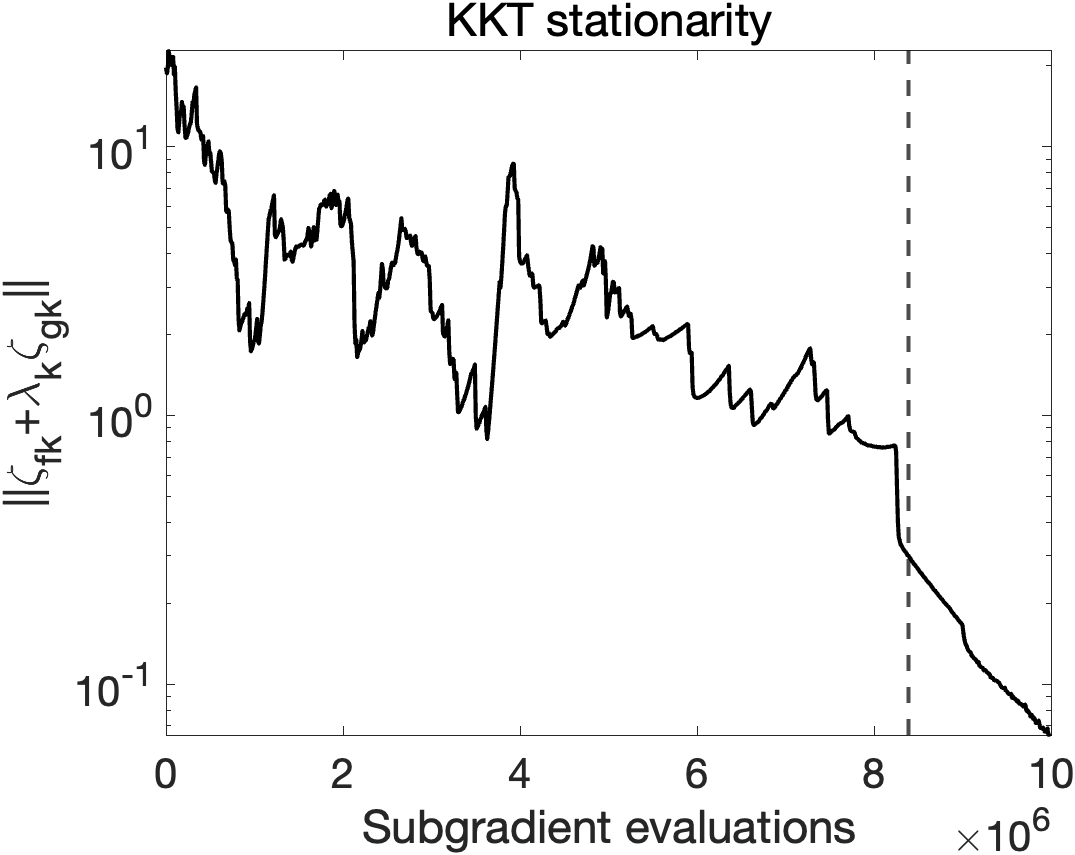}
        \caption{}
        \label{spr2d}
    }
    \end{subfigure}
    \hfill
    \begin{subfigure}{0.3\textwidth}
    {
        \centering
        \includegraphics[width=\textwidth]{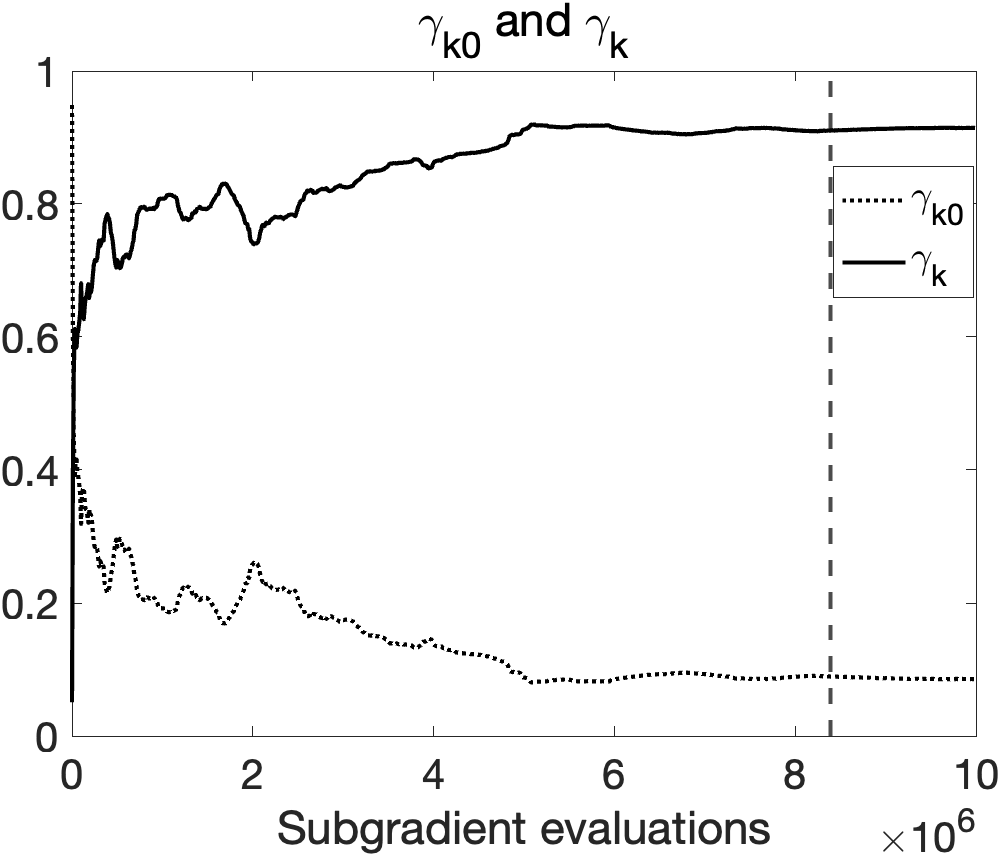}
        \caption{}
        \label{spr2e}
    }
    \end{subfigure}
    \hfill
    \begin{subfigure}{0.3\textwidth}
    {
        \centering
        \includegraphics[width=\textwidth]{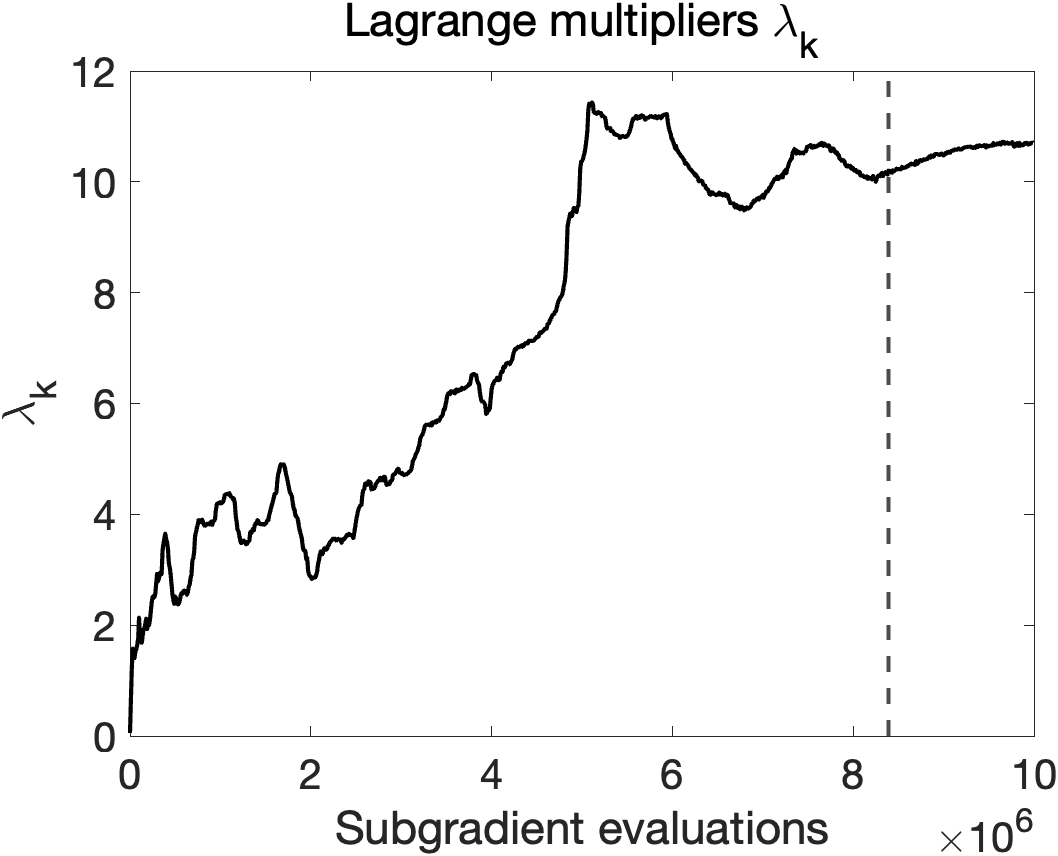}
        \caption{}
        \label{spr2f}
    }
    \end{subfigure}
    \caption{Finding Active KKT Stationarity: $K=10^3,T=10^4,p=91$. Dotted lines show where the stopping criteria applied. $x_{lo}$ is the stationary point near the final iterate.}
    \label{example2}
\end{figure}

\begin{figure}[htbp]
    \centering
    \begin{subfigure}{0.3\textwidth}
    {
        \centering
        \includegraphics[width=\textwidth]{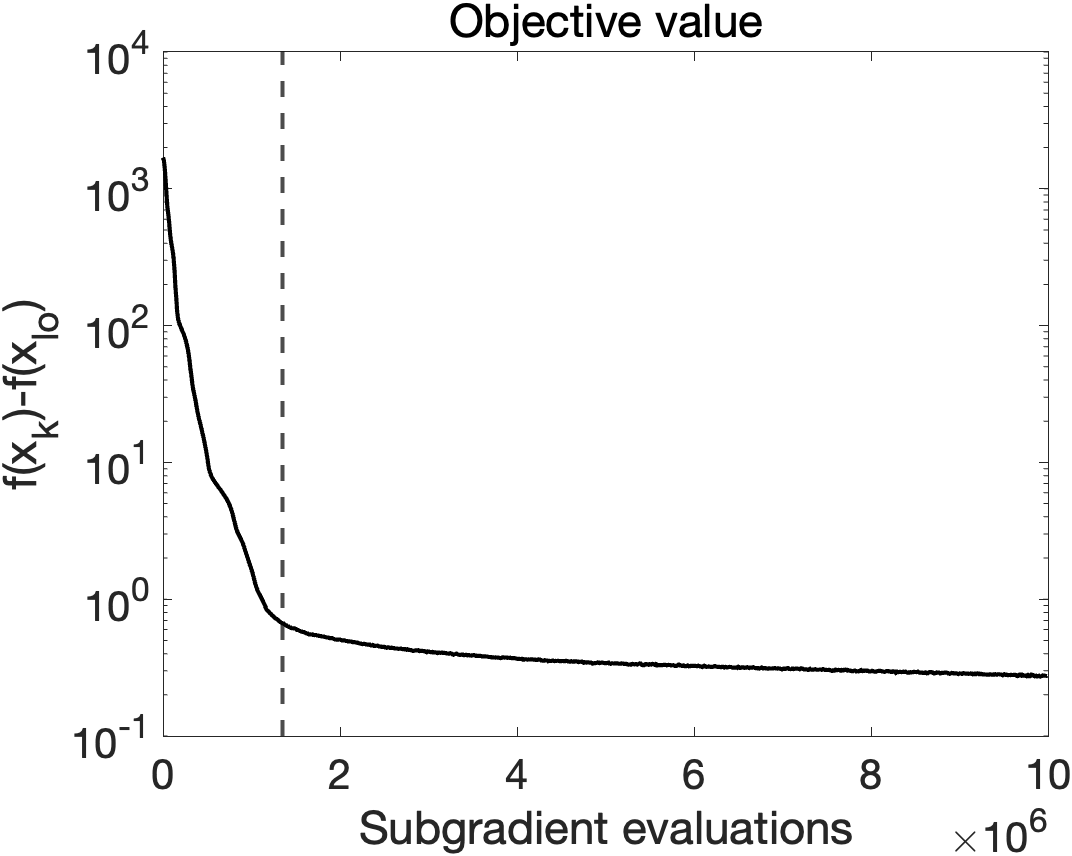}
        \caption{}
        \label{spr3a}
    }
    \end{subfigure}
    \hfill
    \begin{subfigure}{0.3\textwidth}
    {
        \centering
        \includegraphics[width=\textwidth]{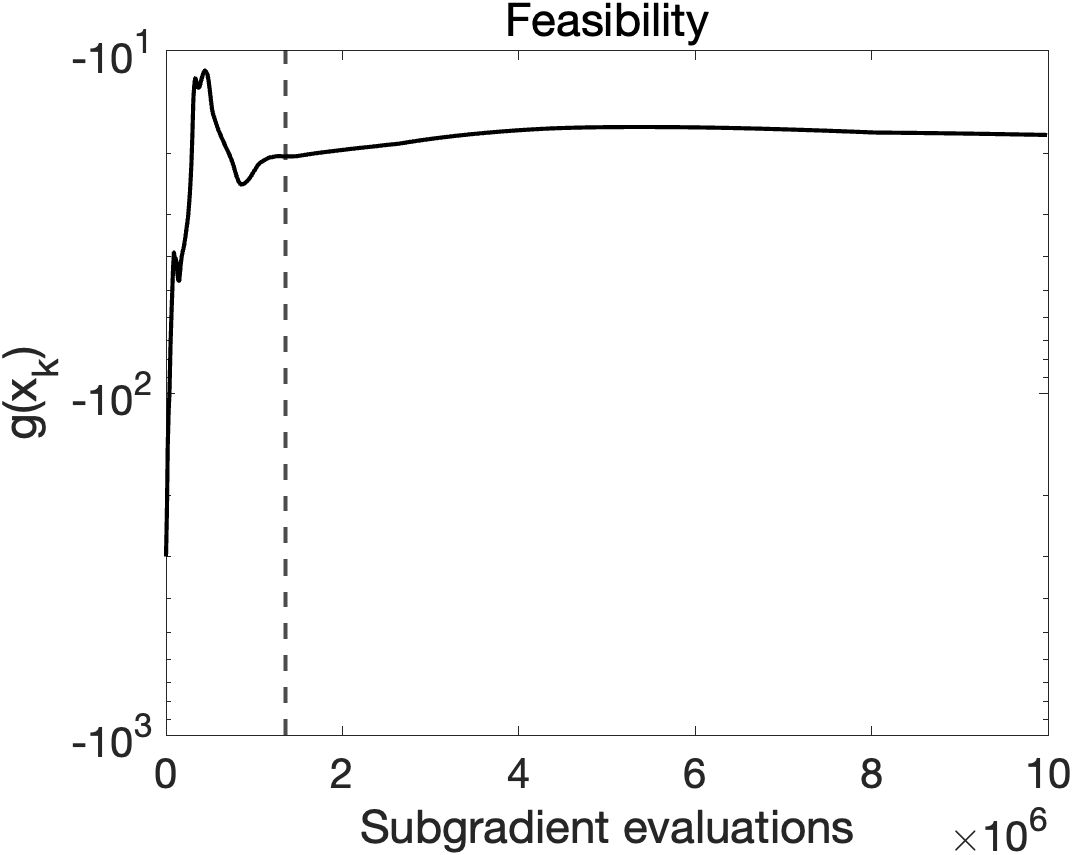}
        \caption{}
        \label{spr3b}
    }
    \end{subfigure}
    \hfill
    \begin{subfigure}{0.3\textwidth}
    {
        \centering
        \includegraphics[width=\textwidth]{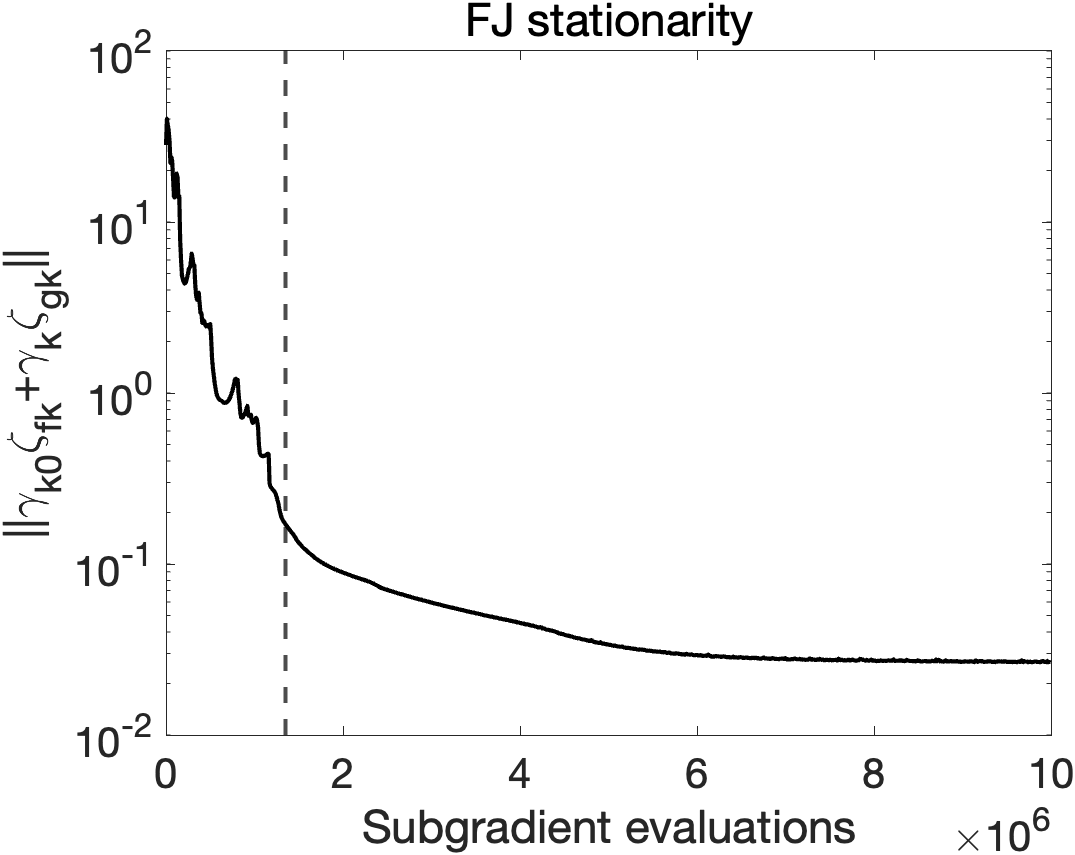}
        \caption{}
        \label{spr3c}
    }
    \end{subfigure}
    \begin{subfigure}{0.3\textwidth}
    {
        \centering
        \includegraphics[width=\textwidth]{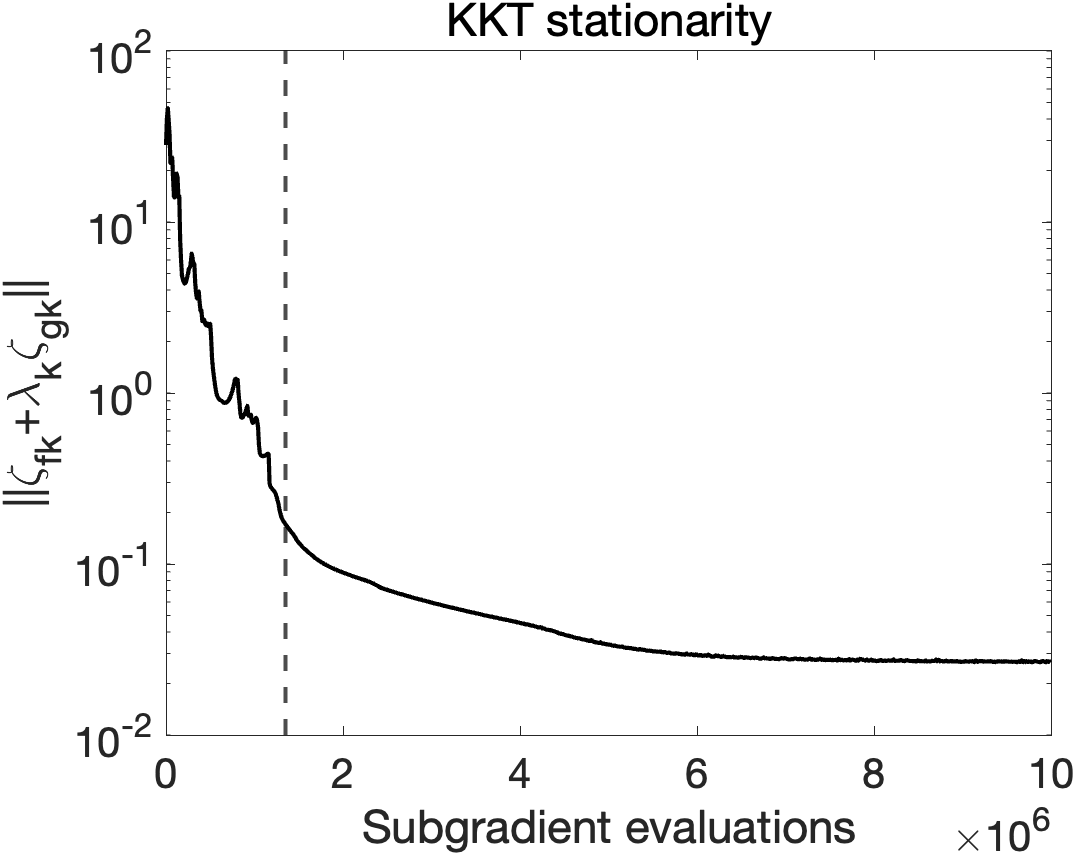}
        \caption{}
        \label{spr3d}
    }
    \end{subfigure}
    \hfill
    \begin{subfigure}{0.3\textwidth}
    {
        \centering
        \includegraphics[width=\textwidth]{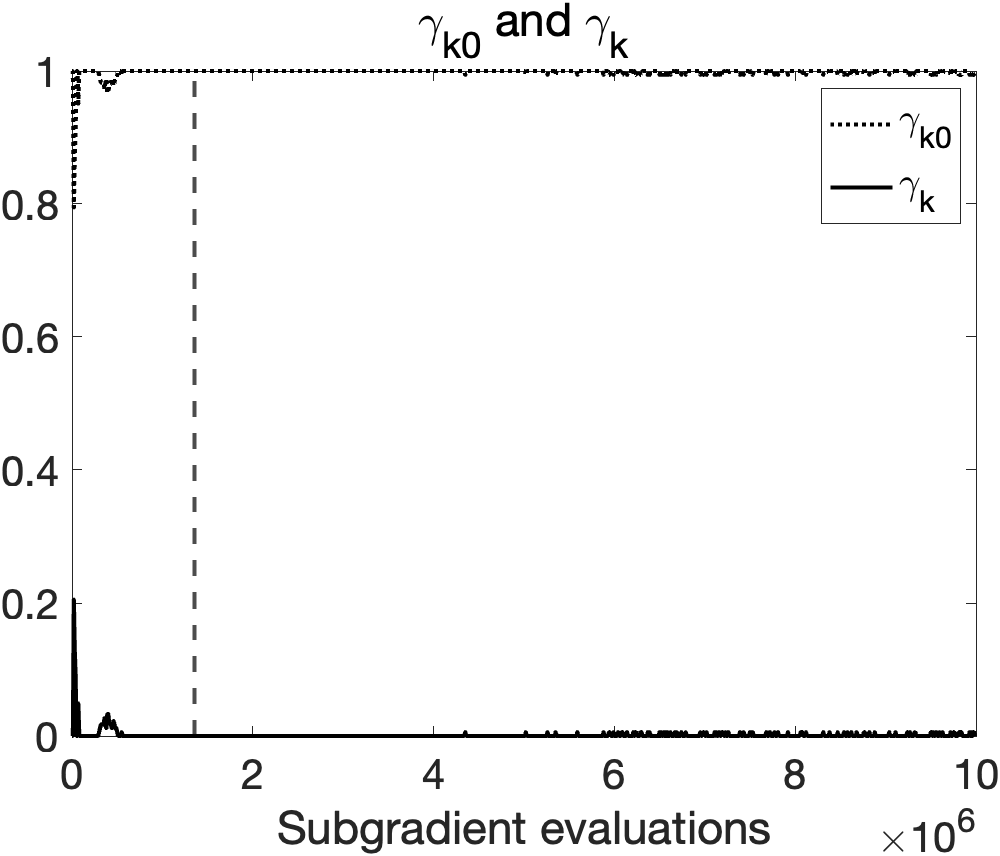}
        \caption{}
        \label{spr3e}
    }
    \end{subfigure}
    \hfill
    \begin{subfigure}{0.3\textwidth}
    {
        \centering
        \includegraphics[width=\textwidth]{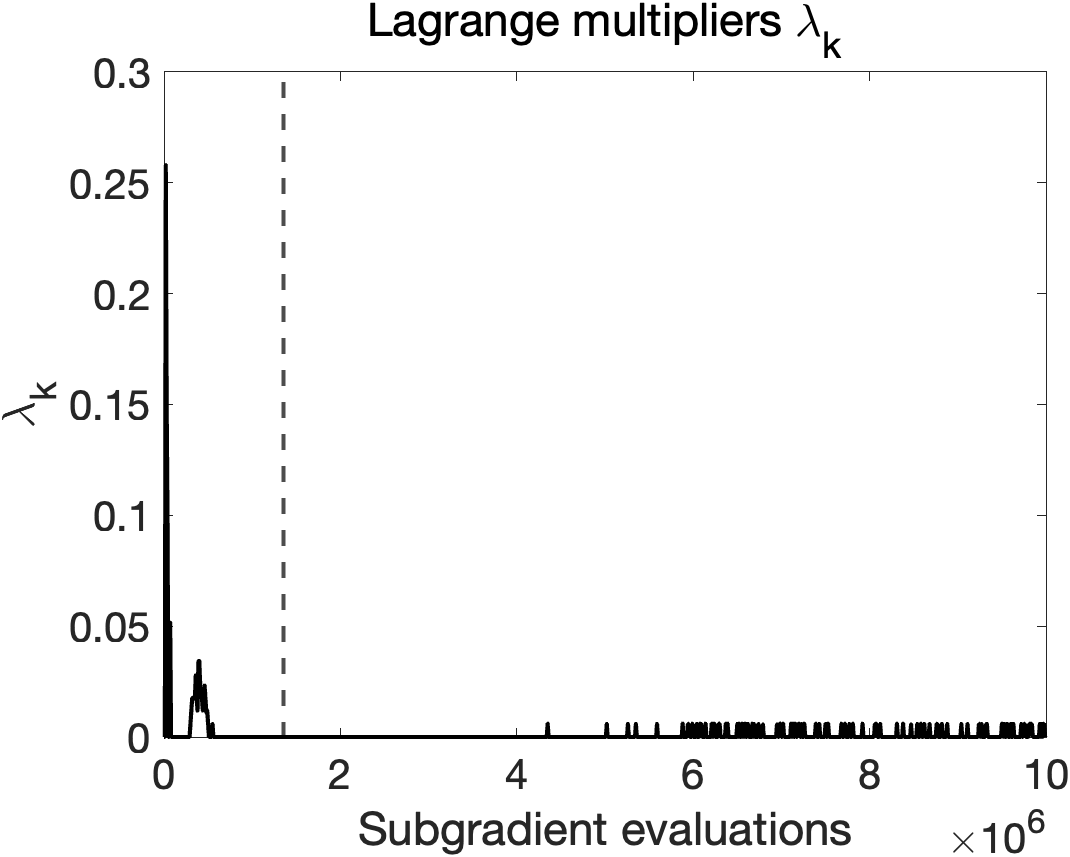}
        \caption{}
        \label{spr3f}
    }
    \end{subfigure}
    \caption{Finding Inactive KKT Stationarity: $K=10^3,T=10^4,p=320$. Dotted lines show where the stopping criteria applied. $x_{lo}$ is the stationary point near the final iterate.}
    \label{example3}
\end{figure}

For each setting of $p$, a sample trajectory of the Proximally Guided Switching Subgradient Method is shown in Figures~\ref{example1}, \ref{example2}, and~\ref{example3}. Statistics on the typical FJ and KKT stationarity levels reached over $50$ trials are provided in Tables~\ref{table1} and~\ref{table2}. Median and variance statistics are included as several experiments (especially those with the potential for MFCQ to fail) had very varied results.

\paragraph{FJ Stationarity}
In the first numeric, we set $p=90$. The example trajectory shown in Figure~\ref{example1} converges to an approximate FJ stationary point of~\eqref{SPR}, which is not an approximate KKT stationary point. Once the stopping criterion is reached, the Lagrangian multiplier estimates diverge rapidly in Figures~\ref{spr1e} and~\ref{spr1f}. Consequently, Figure~\ref{spr1c} shows the FJ stationarity is attained around $10^{-3}$ finally, but Figure~\ref{spr1d} indicates that KKT stationarity is only around $0.5$. Out of the $50$ such trajectories aggregated in Table~\ref{table1}, the shown trajectory is one of three with Lagrange multipliers diverging. This contributes to the larger variance and gap between the mean and median KKT stationarity shown in Table~\ref{table2}.

\paragraph{KKT Stationarity with Active Constraints}

In the second numeric, we set $p=91$. Under this setting, every proximal subproblem satisfies constraint qualification regardless of $x_k$'s location. This is because the subgradient set of $g(x)$ at any $x \in \{x|g(x)=0\}$ contains the zero vector only when every entries of $x$ lies in $(-\infty,2]\cup\{0\}\cup[2,\infty)$, which implies $\sum_i s(x_i)$ is divisible by $3$. As a result, for $p=91$, any $\zeta_g \in g(x)$ taken at the boundary of the constraint set must have size bounded away from zero (ensuring $\sigma$-strong MFCQ). Therefore our inexact proximal point method will always yield an approximate KKT point. We observe this numerically as approximate FJ and KKT stationarity are both reached in Figures~\ref{spr2c} and~\ref{spr2d} and the associated Lagrange multipliers converge to a constant around $10$ in Figures~\ref{spr2e} and~\ref{spr2f}.

\paragraph{KKT Stationarity with Inactive Constraints}
In the third numeric, we set $p=320$. Given this larger value of $p$, we do not expect the constraint to be active or the limit points to be sparse. Complementary slackness at strictly feasible stationary points forces the Lagrange multipliers to equal zero, making FJ and KKT stationarity equivalent. Indeed, Figures~\ref{spr3a} and~\ref{spr3b} show our sample trajectory converges to a strictly feasible local minimum. Figures~\ref{spr3c} and~\ref{spr3d} show that the FJ stationarity and KKT stationarity are equal and converging. As expected, the Lagrange multipliers converge to zero, as shown in Figures~\ref{spr3e} and~\ref{spr3f}.

\begin{table}[htbp]
    \centering\small
    \begin{tabular}{|c|c|c|c|c|c|c|c|}
        \hline
        \multirow{2}*{$KT$} & \multirow{2}*{} & \multicolumn{2}{|c|}{$p=90$} & \multicolumn{2}{|c|}{$p=91$} & \multicolumn{2}{|c|}{$p=320$} \\
        \cline{3-8}
         & & $T=10^3$ & $T=10^4$ & $T=10^3$ & $T=10^4$ & $T=10^3$ & $T=10^4$ \\
        \hline
        \multirow{3}*{$10^5$} & median & $1.174$ & $7.692$ & $1.036$ & $8.327$ & $0.8860$ & $16.11$ \\
        \cline{2-8}
        & mean & $1.239$ & $7.992$ & $1.060$ & $8.272$ & $0.9948$ & $16.34$ \\
        \cline{2-8}
        & var. & $0.4265$ & $1.329$ & $0.1791$ & $1.108$ & $0.2032$ & $11.49$ \\
        \hline
        \multirow{3}*{$10^6$} & median & $0.03230$ & $1.212$ & $0.03370$ & $1.032$ & $0.06256$ & $0.8309$ \\
        \cline{2-8}
        & mean & $0.04563$ & $1.212$ & $0.07497$ & $1.181$ & $0.07423$ & $0.8853$ \\
        \cline{2-8}
        & var. & $2.587$e$-3$ & $0.1937$ & $0.03230$ & $0.3225$ & $1.487$e$-3$ & $0.1453$ \\
        \hline
        \multirow{3}*{$10^7$} & median & $0.03110$ & $2.200$e$-3$ & $0.03140$ & $2.110$e$-3$ & $0.06146$ & $0.01910$ \\
        \cline{2-8}
        & mean & $0.03173$ & $0.04273$ & $0.03130$ & $0.01644$ & $0.06460$ & $0.02941$ \\
        \cline{2-8}
        & var. & $2.545$e$-5$ & $0.01340$ & $2.062$e$-5$ & $3.109$e$-3$ & $2.926$e$-4$& $1.406$e$-3$ \\
        \hline
    \end{tabular}
    \caption{FJ stationarity averaged over $50$ trails.}
    \label{table1}
\end{table}

\begin{table}[htbp]
    \centering\small
    \begin{tabular}{|c|c|c|c|c|c|c|c|}
        \hline
        \multirow{2}*{$KT$} & \multirow{2}*{} & \multicolumn{2}{|c|}{$p=90$} & \multicolumn{2}{|c|}{$p=91$} & \multicolumn{2}{|c|}{$p=320$} \\
        \cline{3-8}
         & & $T=10^3$ & $T=10^4$ & $T=10^3$ & $T=10^4$ & $T=10^3$ & $T=10^4$ \\
        \hline
        \multirow{3}*{$10^5$} & median & $6.904$ & $22.94$ & $6.515$ & $24.44$ & $0.9485$ & $16.39$ \\
        \cline{2-8}
        & mean & $7.810$ & $24.49$ & $6.819$ & $24.94$ & $1.087$ & $16.69$ \\
        \cline{2-8}
        & var. & $20.41$ & $25.92$ & $7.586$ & $23.35$ & $0.2966$ & $15.54$ \\
        \hline
        \multirow{3}*{$10^6$} & median & $0.2812$ & $7.273$ & $0.2987$ & $6.599$ & $0.07857$ & $0.8309$ \\
        \cline{2-8}
        & mean & $0.3855$ & $7.366$ & $0.6012$ & $7.364$ & $0.08914$ & $0.9364$ \\
        \cline{2-8}
        & var. & $0.1502$ & $7.148$ & $1.405$ & $14.31$ & $2.255$e$-3$ & $0.1911$ \\
        \hline
        \multirow{3}*{$10^7$} & median & $0.2740$ & $0.01970$ & $0.2900$ & $0.01752$ & $0.07256$ & $0.0203$ \\
        \cline{2-8}
        & mean & $0.2829$ & $0.5566$ & $0.3146$ & $0.1277$ & $0.07728$ & $0.03169$ \\
        \cline{2-8}
        & var. & $9.349$e$-3$ & $1.777$ & $0.01947$ & $0.1580$ & $4.596$e$-4$ & $1.538$e$-3$ \\
        \hline
    \end{tabular}
    \caption{KKT stationarity averaged over $50$ trails.}
    \label{table2}
\end{table}
    \section{Conclusion and Future Directions}

In this paper, we analyzed an inexact proximal point method using the switching subgradient method as an oracle for nonconvex nonsmooth functional constrained optimization. We derived new convergence rates towards FJ and KKT stationarity while guaranteeing feasibility for our solutions without any reliance on compactness or constraint qualification. The performance of our method for solving sparse phase retrieval problems is consistent with our theoretical expectations. The frequency of constraint qualification failures seen numerically here motivates further works analyzing the performance of nonconvex constrained optimization algorithms both in terms of KKT and FJ convergence.
As additional future directions, stochastic versions of our method could likely be designed and analyzed, like those of \cite{davis2019stochastic,davis2019proximally} from the unconstrained setting or those discussed at the end of Section \ref{switchingsubgradient}. Further, convergence speedups in the presence of structures like local sharpness (see \cite{davis2018subgradient}), strong convexity, or smoothness at the stationary points may be possible.
    
    {\small
    \bibliographystyle{siam}
    \bibliography{bibliography}

\begin{thebibliography}{10}

\bibitem{allgower2012nonlinear}
{\sc F.~Allg{\"o}wer and A.~Zheng}, {\em Nonlinear model predictive control},
  vol.~26, Birkh{\"a}user, 2012.

\bibitem{antoniadis1997wavelets}
{\sc A.~Antoniadis}, {\em Wavelets in statistics: a review}, Journal of the
  Italian Statistical Society, 6 (1997), pp.~97--130.

\bibitem{aravkin2019level}
{\sc A.~Y. Aravkin, J.~V. Burke, D.~Drusvyatskiy, M.~P. Friedlander, and
  S.~Roy}, {\em Level-set methods for convex optimization}, Mathematical
  Programming, 174 (2019), pp.~359--390.

\bibitem{bayandina2018mirror}
{\sc A.~Bayandina, P.~Dvurechensky, A.~Gasnikov, F.~Stonyakin, and A.~Titov},
  {\em Mirror descent and convex optimization problems with non-smooth
  inequality constraints}, in Large-scale and distributed optimization,
  Springer, 2018, pp.~181--213.

\bibitem{bertsekas2015convex}
{\sc D.~Bertsekas}, {\em Convex optimization algorithms}, Athena Scientific,
  2015.

\bibitem{bian2015ipm}
{\sc W.~Bian, X.~Chen, and Y.~Ye}, {\em Complexity analysis of interior point
  algorithms for non-lipschitz and nonconvex minimization}, Mathematical
  Programming,  (2015), pp.~301--327.

\bibitem{birgin2016evaluation}
{\sc E.~G. Birgin, J.~Gardenghi, J.~M. Mart{\'\i}nez, S.~A. Santos, and P.~L.
  Toint}, {\em Evaluation complexity for nonlinear constrained optimization
  using unscaled kkt conditions and high-order models}, SIAM Journal on
  Optimization, 26 (2016), pp.~951--967.

\bibitem{boob2022stochastic}
{\sc D.~Boob, Q.~Deng, and G.~Lan}, {\em Stochastic first-order methods for
  convex and nonconvex functional constrained optimization}, Mathematical
  Programming,  (2022), pp.~1--65.

\bibitem{cartis2011evaluation}
{\sc C.~Cartis, N.~I. Gould, and P.~L. Toint}, {\em On the evaluation
  complexity of composite function minimization with applications to nonconvex
  nonlinear programming}, SIAM Journal on Optimization, 21 (2011),
  pp.~1721--1739.

\bibitem{cartis2013evaluation}
\leavevmode\vrule height 2pt depth -1.6pt width 23pt, {\em On the evaluation
  complexity of cubic regularization methods for potentially rank-deficient
  nonlinear least-squares problems and its relevance to constrained nonlinear
  optimization}, SIAM Journal on Optimization, 23 (2013), pp.~1553--1574.

\bibitem{cartis2014complexity}
\leavevmode\vrule height 2pt depth -1.6pt width 23pt, {\em On the complexity of
  finding first-order critical points in constrained nonlinear optimization},
  Mathematical Programming, 144 (2014), pp.~93--106.

\bibitem{Cartis2015}
\leavevmode\vrule height 2pt depth -1.6pt width 23pt, {\em On the evaluation
  complexity of constrained nonlinear least-squares and general constrained
  nonlinear optimization using second-order methods}, SIAM J. Numer. Anal., 53
  (2015), p.~836–851.

\bibitem{cartis2018second}
{\sc C.~Cartis, N.~I. Gould, and P.~L. Toint}, {\em Second-order optimality and
  beyond: characterization and evaluation complexity in convexly constrained
  nonlinear optimization}, Foundations of Computational Mathematics, 18 (2018),
  pp.~1073--1107.

\bibitem{cartis2020concise}
{\sc C.~Cartis, N.~I. Gould, and P.~L. Toint}, {\em A concise second-order
  complexity analysis for unconstrained optimization using high-order
  regularized models}, Optimization Methods and Software, 35 (2020),
  pp.~243--256.

\bibitem{cartis2020strong}
{\sc C.~Cartis, N.~I. Gould, and P.~L. Toint}, {\em Strong evaluation
  complexity bounds for arbitrary-order optimization of nonconvex nonsmooth
  composite functions}, arXiv preprint arXiv:2001.10802,  (2020).

\bibitem{charisopoulos2021composite}
{\sc V.~Charisopoulos, D.~Davis, M.~D{\'\i}az, and D.~Drusvyatskiy}, {\em
  Composite optimization for robust rank one bilinear sensing}, Information and
  Inference: A Journal of the IMA, 10 (2021), pp.~333--396.

\bibitem{chen2015exact}
{\sc Y.~Chen, Y.~Chi, and A.~J. Goldsmith}, {\em Exact and stable covariance
  estimation from quadratic sampling via convex programming}, IEEE Transactions
  on Information Theory, 61 (2015), pp.~4034--4059.

\bibitem{Clark}
{\sc F.~H. Clarke, Y.~S. Ledyaev, R.~J. Stern, and P.~R. Wolenski}, {\em
  Nonsmooth analysis and control theory}, Springer-Verlag, Berlin, Heidelberg,
  1998.

\bibitem{curtis2012sequential}
{\sc F.~E. Curtis and M.~L. Overton}, {\em A sequential quadratic programming
  algorithm for nonconvex, nonsmooth constrained optimization}, SIAM Journal on
  Optimization, 22 (2012), pp.~474--500.

\bibitem{davis2019stochastic}
{\sc D.~Davis and D.~Drusvyatskiy}, {\em Stochastic model-based minimization of
  weakly convex functions}, SIAM Journal on Optimization, 29 (2019),
  pp.~207--239.

\bibitem{davis2018subgradient}
{\sc D.~Davis, D.~Drusvyatskiy, K.~J. MacPhee, and C.~Paquette}, {\em
  Subgradient methods for sharp weakly convex functions}, Journal of
  Optimization Theory and Applications, 179 (2018), pp.~962--982.

\bibitem{davis2020nonsmooth}
{\sc D.~Davis, D.~Drusvyatskiy, and C.~Paquette}, {\em The nonsmooth landscape
  of phase retrieval}, IMA Journal of Numerical Analysis, 40 (2020),
  pp.~2652--2695.

\bibitem{davis2019proximally}
{\sc D.~Davis and B.~Grimmer}, {\em Proximally guided stochastic subgradient
  method for nonsmooth, nonconvex problems}, SIAM Journal on Optimization, 29
  (2019), pp.~1908--1930.

\bibitem{duchi2019solving}
{\sc J.~C. Duchi and F.~Ruan}, {\em Solving (most) of a set of quadratic
  equalities: Composite optimization for robust phase retrieval}, Information
  and Inference: A Journal of the IMA, 8 (2019), pp.~471--529.

\bibitem{Facchinei2021}
{\sc F.~Facchinei, V.~Kungurtsev, L.~Lampariello, and G.~Scutari}, {\em Ghost
  penalties in nonconvex constrained optimization: Diminishing stepsizes and
  iteration complexity}, Mathematics of Operations Research, 46 (2021),
  pp.~595--627.

\bibitem{fan2009network}
{\sc J.~Fan, Y.~Feng, and Y.~Wu}, {\em Network exploration via the adaptive
  lasso and scad penalties}, The annals of applied statistics, 3 (2009),
  p.~521.

\bibitem{fan2001variable}
{\sc J.~Fan and R.~Li}, {\em Variable selection via nonconcave penalized
  likelihood and its oracle properties}, Journal of the American statistical
  Association, 96 (2001), pp.~1348--1360.

\bibitem{frostig2015regularizing}
{\sc R.~Frostig, R.~Ge, S.~Kakade, and A.~Sidford}, {\em Un-regularizing:
  approximate proximal point and faster stochastic algorithms for empirical
  risk minimization}, in International Conference on Machine Learning, PMLR,
  2015, pp.~2540--2548.

\bibitem{gasso2009recovering}
{\sc G.~Gasso, A.~Rakotomamonjy, and S.~Canu}, {\em Recovering sparse signals
  with a certain family of nonconvex penalties and dc programming}, IEEE
  Transactions on Signal Processing, 57 (2009), pp.~4686--4698.

\bibitem{gower2019sgd}
{\sc R.~M. Gower, N.~Loizou, X.~Qian, A.~Sailanbayev, E.~Shulgin, and
  P.~Richt{\'a}rik}, {\em Sgd: General analysis and improved rates}, in
  International Conference on Machine Learning, PMLR, 2019, pp.~5200--5209.

\bibitem{grapiglia2021complexity}
{\sc G.~N. Grapiglia and Y.-x. Yuan}, {\em On the complexity of an augmented
  lagrangian method for nonconvex optimization}, IMA Journal of Numerical
  Analysis, 41 (2021), pp.~1546--1568.

\bibitem{grimmer2019convergence}
{\sc B.~Grimmer}, {\em Convergence rates for deterministic and stochastic
  subgradient methods without lipschitz continuity}, SIAM Journal on
  Optimization, 29 (2019), pp.~1350--1365.

\bibitem{radial1}
\leavevmode\vrule height 2pt depth -1.6pt width 23pt, {\em Radial duality part
  i: foundations}, Mathematical Programming,  (2021).

\bibitem{radial2}
\leavevmode\vrule height 2pt depth -1.6pt width 23pt, {\em Radial duality part
  ii: applications and algorithms}, Mathematical Programming,  (2021).

\bibitem{guler1992new}
{\sc O.~G{\"u}ler}, {\em New proximal point algorithms for convex
  minimization}, SIAM Journal on Optimization, 2 (1992), pp.~649--664.

\bibitem{hare2009computing}
{\sc W.~Hare and C.~Sagastiz{\'a}bal}, {\em Computing proximal points of
  nonconvex functions}, Mathematical Programming, 116 (2009), pp.~221--258.

\bibitem{hare2010redistributed}
\leavevmode\vrule height 2pt depth -1.6pt width 23pt, {\em A redistributed
  proximal bundle method for nonconvex optimization}, SIAM Journal on
  Optimization, 20 (2010), pp.~2442--2473.

\bibitem{hinder2018worst}
{\sc O.~Hinder and Y.~Ye}, {\em Worst-case iteration bounds for log barrier
  methods for problems with nonconvex constraints}, arXiv preprint
  arXiv:1807.00404,  (2018).

\bibitem{john2014extremum}
{\sc F.~John}, {\em Extremum problems with inequalities as subsidiary
  conditions}, in Traces and emergence of nonlinear programming, Springer,
  2014, pp.~197--215.

\bibitem{jourani1994constraint}
{\sc A.~Jourani}, {\em Constraint qualifications and lagrange multipliers in
  nondifferentiable programming problems}, Journal of Optimization Theory and
  Applications, 81 (1994), pp.~533--548.

\bibitem{kim2008smoothly}
{\sc Y.~Kim, H.~Choi, and H.-S. Oh}, {\em Smoothly clipped absolute deviation
  on high dimensions}, Journal of the American Statistical Association, 103
  (2008), pp.~1665--1673.

\bibitem{kong2019complexity}
{\sc W.~Kong, J.~G. Melo, and R.~D. Monteiro}, {\em Complexity of a quadratic
  penalty accelerated inexact proximal point method for solving linearly
  constrained nonconvex composite programs}, SIAM Journal on Optimization, 29
  (2019), pp.~2566--2593.

\bibitem{lacostejulien2012simpler}
{\sc S.~Lacoste-Julien, M.~Schmidt, and F.~Bach}, {\em A simpler approach to
  obtaining an o(1/t) convergence rate for the projected stochastic subgradient
  method}, 2012.

\bibitem{lan2019accelerated}
{\sc G.~Lan and Y.~Yang}, {\em Accelerated stochastic algorithms for nonconvex
  finite-sum and multiblock optimization}, SIAM Journal on Optimization, 29
  (2019), pp.~2753--2784.

\bibitem{lan2020algorithms}
{\sc G.~Lan and Z.~Zhou}, {\em Algorithms for stochastic optimization with
  function or expectation constraints}, Computational Optimization and
  Applications, 76 (2020), pp.~461--498.

\bibitem{lin2018level}
{\sc Q.~Lin, S.~Nadarajah, and N.~Soheili}, {\em A level-set method for convex
  optimization with a feasible solution path}, SIAM Journal on Optimization, 28
  (2018), pp.~3290--3311.

\bibitem{ma2020quadratically}
{\sc R.~Ma, Q.~Lin, and T.~Yang}, {\em Quadratically regularized subgradient
  methods for weakly convex optimization with weakly convex constraints}, in
  International Conference on Machine Learning, PMLR, 2020, pp.~6554--6564.

\bibitem{MANGASARIAN196737}
{\sc O.~Mangasarian and S.~Fromovitz}, {\em The fritz john necessary optimality
  conditions in the presence of equality and inequality constraints}, Journal
  of Mathematical Analysis and Applications, 17 (1967), pp.~37--47.

\bibitem{mangasarian1967fritz}
{\sc O.~L. Mangasarian and S.~Fromovitz}, {\em The fritz john necessary
  optimality conditions in the presence of equality and inequality
  constraints}, Journal of Mathematical Analysis and applications, 17 (1967),
  pp.~37--47.

\bibitem{martinet1970regularisation}
{\sc B.~Martinet}, {\em Regularisation d'inequations variationelles par
  approximations successives}, Revue Francaise d'informatique et de Recherche
  operationelle, 4 (1970), pp.~154--159.

\bibitem{moreau1965proximite}
{\sc J.-J. Moreau}, {\em Proximit{\'e} et dualit{\'e} dans un espace
  hilbertien}, Bulletin de la Soci{\'e}t{\'e} math{\'e}matique de France, 93
  (1965), pp.~273--299.

\bibitem{Nocedal2006}
{\sc J.~Nocedal and S.~J. Wright}, {\em Sequential Quadratic Programming},
  Springer New York, New York, NY, 2006, pp.~529--562.

\bibitem{nouiehed2018convergence}
{\sc M.~Nouiehed, J.~D. Lee, and M.~Razaviyayn}, {\em Convergence to
  second-order stationarity for constrained non-convex optimization}, arXiv
  preprint arXiv:1810.02024,  (2018).

\bibitem{paquette2018catalyst}
{\sc C.~Paquette, H.~Lin, D.~Drusvyatskiy, J.~Mairal, and Z.~Harchaoui}, {\em
  Catalyst for gradient-based nonconvex optimization}, in International
  Conference on Artificial Intelligence and Statistics, PMLR, 2018,
  pp.~613--622.

\bibitem{peng2020computation}
{\sc D.~Peng and X.~Chen}, {\em Computation of second-order directional
  stationary points for group sparse optimization}, Optimization Methods and
  Software, 35 (2020), pp.~348--376.

\bibitem{pieper2022nonconvex}
{\sc K.~Pieper and A.~Petrosyan}, {\em Nonconvex regularization for sparse
  neural networks}, Applied and Computational Harmonic Analysis,  (2022).

\bibitem{articlepolyak}
{\sc B.~Polyak}, {\em A general method for solving extremum problems}, Soviet
  Mathematics. Doklady, 8 (1967).

\bibitem{salzo2012inexact}
{\sc S.~Salzo and S.~Villa}, {\em Inexact and accelerated proximal point
  algorithms}, Journal of Convex analysis, 19 (2012), pp.~1167--1192.

\bibitem{tian2021neyman}
{\sc Y.~Tian and Y.~Feng}, {\em Neyman-pearson multi-class classification via
  cost-sensitive learning}, arXiv preprint arXiv:2111.04597,  (2021).

\bibitem{Wang2017}
{\sc X.~Wang, S.~Ma, and Y.-X. Yuan}, {\em Penalty methods with stochastic
  approximation for stochastic nonlinear programming}, Mathematics of
  Computation, 86 (2017), pp.~1793--1820.

\bibitem{wen2018survey}
{\sc F.~Wen, L.~Chu, P.~Liu, and R.~C. Qiu}, {\em A survey on nonconvex
  regularization-based sparse and low-rank recovery in signal processing,
  statistics, and machine learning}, IEEE Access, 6 (2018), pp.~69883--69906.

\bibitem{wen2017efficient}
{\sc F.~Wen, L.~Pei, Y.~Yang, W.~Yu, and P.~Liu}, {\em Efficient and robust
  recovery of sparse signal and image using generalized nonconvex
  regularization}, IEEE Transactions on Computational Imaging, 3 (2017),
  pp.~566--579.

\bibitem{weston1998multi}
{\sc J.~Weston and C.~Watkins}, {\em Multi-class support vector machines},
  tech. rep., Citeseer, 1998.

\bibitem{xie2009scad}
{\sc H.~Xie and J.~Huang}, {\em Scad-penalized regression in high-dimensional
  partially linear models}, The Annals of Statistics, 37 (2009), pp.~673--696.

\bibitem{yu2017design}
{\sc Z.~Yu, P.~Cui, and J.~L. Crassidis}, {\em Design and optimization of
  navigation and guidance techniques for mars pinpoint landing: Review and
  prospect}, Progress in Aerospace Sciences, 94 (2017), pp.~82--94.

\bibitem{zeng2022moreau}
{\sc J.~Zeng, W.~Yin, and D.-X. Zhou}, {\em Moreau envelope augmented
  lagrangian method for nonconvex optimization with linear constraints},
  Journal of Scientific Computing, 91 (2022), pp.~1--36.

\bibitem{zhang2020meta}
{\sc H.~Zhang, S.-J. Li, H.~Zhang, Z.-Y. Yang, Y.-Q. Ren, L.-Y. Xia, and
  Y.~Liang}, {\em Meta-analysis based on nonconvex regularization}, Scientific
  reports, 10 (2020), pp.~1--16.

\end{thebibliography}
    
    \appendix
    \section{Deferred Proofs}

\subsection{Symmetric KKT Proof of Theorem~\ref{theorem_mainKKT}}
According to Lemma \ref{lemmafea}, our iterates are always feasible, that is $g(x_k) \leq 0$, for the main problem \eqref{mainproblem} before $\hat x_{k+1}$ is an $\epsilon$-KKT point (which will imply $x_{k}$ is an $(\epsilon,\epsilon)$-KKT point). For each iteration $k$, let $\lambda_k$ denote the optimal Lagrange multiplier in \eqref{KKToriginal} for the proximal subproblem \eqref{subproblem}. We denote the Lagrange function for each subproblem~\eqref{subproblem} as
\begin{equation}
    L_k(x)=F_k(x)+\lambda_kG_k(x)=f(x)+\frac{\hat{\rho}}{2}\|x-x_k\|^2+\lambda_k\left(g(x)+\frac{\hat{\rho}}{2}\|x-x_k\|^2\right)\ .
    \label{lagrange2}
\end{equation}
Without loss of generality, suppose $\lambda_k \geq 0$. According to KKT conditions \eqref{KKToriginal}, there exists $\hat{\zeta}_{Fk} \in \partial F_k(\hat{x}_{k+1})$ and $\hat{\zeta}_{Gk} \in \partial G_k(\hat{x}_{k+1})$ which satisfies
\begin{equation}
    \hat{\zeta}_{Fk}+\lambda_k \hat{\zeta}_{Gk} \in -N_X(\hat{x}_{k+1})\ .
    \label{KKTcondition1}
\end{equation}
Since $L_k(x)$ is $(1+\lambda_k)(\hat{\rho}-\rho)$-strongly convex, we have
\begin{equation}
\begin{aligned}
    F_k(x_k)+\lambda_kG_k(x_k) \geq &F_k(\hat{x}_{k+1})+\lambda_kG_k(\hat{x}_{k+1})+(\hat{\zeta}_{Fk}+\lambda_k \hat{\zeta}_{Gk})^T(x_k-\hat{x}_{k+1}) \\
    &+\frac{(1+\lambda_k)(\hat{\rho}-\rho)}{2}\|\hat{x}_{k+1}-x_k\|^2\ .\nonumber
\end{aligned}
\end{equation}
According to KKT conditions \eqref{KKToriginal}, we also have $\lambda_kG_k(\hat{x}_{k+1})=0$. By \eqref{KKTcondition1} and since $x_k \in X$, we know $(\hat{\zeta}_{Fk}+\lambda_k \hat{\zeta}_{Gk})^T(x_k-\hat{x}_{k+1}) \geq 0$. Since $g(x_k) \leq 0$ from Lemma \ref{lemmafea}, the previous inequality becomes
\begin{equation}
    f(x_k) \geq F_k(\hat{x}_{k+1})+\frac{\hat{\rho}-\rho}{2}\|\hat{x}_{k+1}-x_k\|^2\ .\nonumber
\end{equation}
Since $x_{k+1}$ is a $(\tau, \tau)$-solution for problem \eqref{subproblem}, $F_k(x_{k+1})-F_k(\hat{x}_{k+1}) \leq \tau$, then
\begin{align}
    f(x_k) &\geq \left(f(x_{k+1})+\frac{\hat{\rho}}{2}\|x_{k+1}-x_k\|^2-\tau\right)+\frac{(1+\lambda_k)(\hat{\rho}-\rho)}{2}\|\hat{x}_{k+1}-x_k\|^2 \nonumber\\
    &\geq f(x_{k+1})-\tau+\frac{(1+\lambda_k)(\hat{\rho}-\rho)}{2}\|\hat{x}_{k+1}-x_k\|^2\ .\nonumber
\end{align}
Thus we attain a lower bound for the descent of each step as
\begin{equation}
    f(x_k)-f(x_{k+1}) \geq \frac{(1+\lambda_k)(\hat{\rho}-\rho)}{2}\|\hat{x}_{k+1}-x_k\|^2-\tau\ .
    \label{flborigin2}
\end{equation}
Let $\hat{\zeta}_{fk}=\hat{\zeta}_{Fk}-\hat{\rho}(\hat{x}_{k+1}-x_k) \in \partial f(\hat{x}_{k+1})$, $\hat{\zeta}_{gk}=\hat{\zeta}_{Gk}-\hat{\rho}(\hat{x}_{k+1}-x_k) \in \partial g(\hat{x}_{k+1})$. According to Lemma \ref{lemmastationarity}, before $\hat x_{k+1}$ is an $\epsilon$-KKT point, $\|\hat{x}_{k+1}-x_k\|>\frac{\epsilon}{(1+\lambda_k)\hat{\rho}}$, then our choice of $\tau$ ensures
\begin{align}
    f(x_k)-f(x_{k+1}) &\geq \frac{(1+\lambda_k)(\hat{\rho}-\rho)}{2}\|\hat{x}_{k+1}-x_k\|^2-\tau \nonumber\\
    &>\frac{(1+\lambda_k)(\hat{\rho}-\rho)}{2} \frac{\epsilon^2}{(1+\lambda_k)^2\hat{\rho}^2}-\frac{(\hat{\rho}-\rho)\epsilon^2}{4(1+B)^2\hat{\rho}(2\hat{\rho}-\rho)} \nonumber\\
    &>\frac{(\hat{\rho}-\rho)\epsilon^2}{2(1+\lambda_k)\hat{\rho}^2}-\frac{(\hat{\rho}-\rho)\epsilon^2}{4(1+B)\hat{\rho}(2\hat{\rho}-\rho)} \nonumber\\
    &\geq \frac{(\hat{\rho}-\rho)\epsilon^2}{2(1+B)\hat{\rho}^2}-\frac{(\hat{\rho}-\rho)\epsilon^2}{4(1+B)\hat{\rho}^2} \nonumber\\
    &=\frac{(\hat{\rho}-\rho)\epsilon^2}{4(1+B)\hat{\rho}^2}\ . \nonumber
\end{align}
By Assumption \ref{assumption2}, we could give an upper bound for the number of total iterations $K$ as
\begin{equation}
    K<\frac{4(1+B)\hat{\rho}^2(f(x_0)-f_{lb})}{(\hat{\rho}-\rho)\epsilon^2}\ .\nonumber
\end{equation}
Consequently, Algorithm~\ref{algorithm_main} (which uses Algorithm~\ref{algorithm_sub} for $T$ steps as an oracle each iteration) will identify an $(\epsilon,\epsilon)$-KKT point using at most $KT=O(1/\epsilon^4)$ total subgradient evaluations of either the objective or constraints.

\subsection{Symmetric KKT Case of Lemma~\ref{lemmastationarity}'s Proof}
Given constraint qualification, necessarily, the KKT conditions are satisfied for the proximal subproblem~\eqref{subproblem} by some $\lambda_k \geq 0$, $\hat{\zeta}_{Fk} \in \partial F_k(\hat{x}_{k+1})$ and $\hat{\zeta}_{Gk} \in \partial G_k(\hat{x}_{k+1})$. By the sum rule, let $\hat{\zeta}_{fk}=\hat{\zeta}_{Fk}-\hat{\rho}(\hat{x}_{k+1}-x_k) \in \partial f(\hat{x}_{k+1})$ and $\hat{\zeta}_{gk}=\hat{\zeta}_{Gk}-\hat{\rho}(\hat{x}_{k+1}-x_k) \in \partial g(\hat{x}_{k+1})$. The KKT conditions for the proximal subproblem ensure some $\nu \in N_X(\hat{x}_{k+1})$ has
\begin{equation}
    \hat{\zeta}_{fk}+\hat{\rho}(\hat{x}_{k+1}-x_k)+\lambda_k\left[\hat{\zeta}_{gk}+\hat{\rho}(\hat{x}_{k+1}-x_k)\right]=-\nu \ . \nonumber
\end{equation}
When $\|\hat{x}_{k+1}-x_k\| \leq \frac{\epsilon}{\hat{\rho}(1+\lambda_k)}$, it follows that $\|\hat{\zeta}_{fk}+\lambda_k\hat{\zeta}_{gk}+\nu\|=\hat{\rho}(1+\lambda_k)\|\hat{x}_{k+1}-x_k\| \leq \epsilon$, establishing the first approximate KKT condition~\eqref{KKT1_1}. We verify the second approximate KKT condition~\eqref{KKT1_2} in two cases: When $\lambda_k=0$, trivially $|\lambda_k g(\hat{x}_{k+1})|=0$. When $\lambda_k$ is positive, the KKT conditions for the proximal subproblem require $G_k(\hat{x}_{k+1})=0$.
Hence $0 \geq g(\hat{x}_{k+1})=-\frac{\hat{\rho}}{2}\|\hat{x}_{k+1}-x_k\|^2 \geq -\frac{\epsilon^2}{2\hat{\rho}(1+\lambda_k)^2}$\ .
Therefore
\begin{equation}
    |\lambda_kg(\hat{x}_{k+1})| \leq \frac{\epsilon^2\lambda_k}{2\hat{\rho}(1+\lambda_k)^2} \leq \frac{\epsilon^2}{2\hat{\rho}}<\epsilon^2 \ .\nonumber
\end{equation}

\subsection{Symmetric KKT Case of Lemma~\ref{lemmafea}'s Proof}
Assume $G_k(x_k)=g(x_k) \leq 0$. Necessarily, there exists an optimal dual variable for the subproblem $\lambda_k$. For the Lagrange function $L_k(x) = F_k(x) + \lambda_k G_k(x)$, which $\hat x_{k+1}$ minimizes, its $(1+\lambda_k)(\hat{\rho}-\rho)$-strong convexity ensures
\begin{align*}
    F_k(x_{k+1})+\lambda_kG_k(x_{k+1})
    \geq &F_k(\hat{x}_{k+1})+\lambda_kG_k(\hat{x}_{k+1})+(\hat{\zeta}_{Fk}+\lambda_k \hat{\zeta}_{Gk})^T(x_{k+1}-\hat{x}_{k+1}) \\
    &+\frac{(1+\lambda_k)(\hat{\rho}-\rho)}{2}\|x_{k+1}-\hat{x}_{k+1}\|^2\ .
\end{align*}
Note the KKT conditions ensure that $\lambda_kG_k(\hat{x}_{k+1})=0$ 
 and that $\hat{\zeta}_{Fk}+\lambda_k \hat{\zeta}_{Gk} \in -N_X(\hat{x}_{k+1})$, from which one can conclude $x_{k+1} \in X$ must have $(\zeta_{Fk}+\lambda_k \zeta_{Gk})^T(x_{k+1}-\hat{x}_{k+1}) \geq 0$. Then the above inequality simplifies to
\begin{align*}
    F_k(x_{k+1}) - F_k(\hat{x}_{k+1}) +\lambda_kG_k(x_{k+1})
    \geq \frac{(1+\lambda_k)(\hat{\rho}-\rho)}{2}\|x_{k+1}-\hat{x}_{k+1}\|^2\ .
\end{align*}
The proposed value of $T$ ensures $x_{k+1}$ is a $(\tau,\tau)$-optimal solution for the subproblem \eqref{subproblem}, i.e., $F_k(x_{k+1})-F_k(\hat{x}_{k+1}) \leq \tau$ and $G_k(x_{k+1}) \leq \tau$. Hence
\begin{align}
    \|\hat{x}_{k+1}-x_{k+1}\| &\leq \sqrt{\frac{2\tau}{\hat{\rho}-\rho}}\ .\nonumber
\end{align}
Assume $\hat x_{k+1}$ is not an $\epsilon$-KKT point, Lemma~\ref{lemmastationarity} implies $\|\hat{x}_{k+1}-x_k\|>\frac{\epsilon}{(1+B)\hat{\rho}}$. Thus

\begin{equation}
    \|x_{k+1}-x_k\|^2 \geq \frac{1}{2}\|\hat{x}_{k+1}-x_k\|^2-\|\hat{x}_{k+1}-x_{k+1}\|^2>\frac{\epsilon^2}{2(1+B)^2\hat{\rho}^2}-\frac{2\tau}{\hat{\rho}-\rho} \ .\nonumber
\end{equation}
By our selection of $\tau=\frac{(\hat{\rho}-\rho)\epsilon^2}{4(1+B)^2\hat{\rho}(2\hat{\rho}-\rho)}$ as in \eqref{paraKKT}, every iteration prior to finding an  $\epsilon$-KKT point must have
\begin{equation}
    \|x_{k+1}-x_k\|^2>\frac{(\hat{\rho}-\rho)\epsilon^2}{2(1+B)^2\hat{\rho}^2(2\hat{\rho}-\rho)}\ .
    \label{stoppingcriteriaKKT2}
\end{equation}
Therefore $g(x_{k+1}) \leq 0$ is inductively ensured if $g(x_k) \leq 0$ and $\hat x_{k+1}$ is not an $\epsilon$-KKT point as
\begin{align*}
    g(x_{k+1})=G(x_{k+1})-\frac{\hat{\rho}}{2}\|x_{k+1}-x_k\|^2
    \leq \tau-\frac{\hat{\rho}}{2}\frac{(\hat{\rho}-\rho)\epsilon^2}{2(1+B)^2\hat{\rho}^2(2\hat{\rho}-\rho)}=0 \ .
\end{align*}

{\color{black}
\section{Efficient Projection on Proximal SCAD Subproblems} \label{app:SCAD_proj}
For the special case considered in Section~\ref{subsec:vignette} where $g(x) = \sum_{i=1}^n s(x_i) - p$, orthogonal projections for each proximal subproblem can be done efficiently. For any proximal center $\bar x$ with $g(\bar x)<0$, the convex feasible region of the proximal subproblem is given by
$$ Q = \left\{ x \Bigg| g(x) + \frac{1}{2\alpha}\|x-\bar x\|^2\right\} = \left\{ x \Bigg| \sum_{i=1}^n \left[s(x_i) + \frac{1}{2\alpha}(x_i-\bar x_i)^2\right] \leq p\right\} \ . $$
Then the orthogonal projection of any $z\in\mathbb{R}^n$ onto $Q$ can be reformulated as
\begin{align*}
    \begin{cases}
        \min_{x\in\mathbb{R}^n} & \|x - z\|^2\\
        \mathrm{s.t.} & \sum_{i=1}^n \left[s(x_i) + \frac{1}{2\alpha}(x_i-\bar x_i)^2\right] \leq p
    \end{cases} &= \min_x \max_{\lambda\geq 0} \|x - z\|^2 + \lambda\left\{\sum_{i=1}^n \left[s(x_i) + \frac{1}{2\alpha}(x_i-\bar x_i)^2\right] - p\right\}\\
    &= \max_{\lambda\geq 0} \min_x \|x - z\|^2 + \lambda\left\{\sum_{i=1}^n \left[s(x_i) + \frac{1}{2\alpha}(x_i-\bar x_i)^2\right] - p\right\}\\
    &= \max_{\lambda\geq 0} \sum_{i=1}^n \min_{x_i} \left\{(x_i-z_i)^2+ \lambda\left[s(x_i) + \frac{1}{2\alpha}(x_i-\bar x_i)^2\right] - p\right\}
\end{align*}
where the first equality is the Lagrange primal formulation, the second uses strong duality as $\bar x$ is a strictly feasible Slater point, and the third reorganizes terms to make the inner minimization separable over the coordinates of $x_i$. For a given scalar $\lambda\geq 0$, this last inner minimization problem is a proximal step on the SCAD function, which can be cheaply solved in closed form for each $x_i$. This proximal evaluation will often yield $x_i=0$, giving sparsity overall. Consequently, orthogonal projection onto proximal SCAD constraints reduces to one-dimensional maximization over $\lambda$, which can be done rapidly via any root-finding procedure.
}
    }
\end{document}